\documentclass[11pt]{article}

\usepackage{authblk}
\usepackage{url}
\usepackage[margin=1in,footskip=0.25in]{geometry}
\usepackage[labelfont=bf]{caption}
\usepackage{graphicx}
\usepackage{amsmath,amssymb,amsfonts,amsmath}
\usepackage{graphicx}
\usepackage{subfigure}
\usepackage{caption}
\usepackage{tikz}
\usepackage{float}
\usepackage[dvipdf]{epsfig}
\usepackage{float}
\usepackage{lipsum}
\usepackage{wrapfig}
\usepackage{cite}
\usepackage{color}
\usepackage{easylist}
\usepackage{lineno}
\usepackage{booktabs}
\usepackage{hyperref}
\usepackage{cleveref}
\usepackage{multirow}
\usepackage{hhline}
\usepackage{soul}
\usepackage{blindtext}
\usepackage{threeparttable}
\usepackage[bbgreekl]{mathbbol}

\usepackage{amsmath,amssymb,amsfonts,amsmath}
\usepackage{graphicx}
\usepackage{lineno,hyperref}
\usepackage{subfigure}
\usepackage{caption}
\usepackage{float}
\usepackage{mathabx}
\usepackage{tikz}
\usepackage{algorithmicx}
\usepackage{algpseudocode}
\usepackage[dvipdf]{epsfig}

\usepackage{graphicx}

\usepackage{dcolumn}
\usepackage{bm}

\usepackage[utf8]{inputenc}
\usepackage[T1]{fontenc}
\usepackage{mathptmx}
\usepackage{etoolbox}
\usepackage{subfigure}
\usepackage{hyperref}
\usepackage{color}
\usepackage{amsmath,amssymb,amsfonts,amsmath}
\usepackage{textcase}

\newcommand{\eqn}[1]{
\begin{eqnarray}
#1
\end{eqnarray}
}

\newcommand{\norm}[1]{
\left\lVert #1 \right\rVert
}

\newcommand{\ip}[1]{
\left\langle #1 \right\rangle
}

\newcommand{\bs}[1]{
\boldsymbol #1
}

\newcommand{\pz}{
\partial_{z}
}

\makeatletter
\def\@email#1#2{%
\endgroup
\patchcmd{\titleblock@produce}
{\frontmatter@RRAPformat}
{\frontmatter@RRAPformat{\produce@RRAP{*#1\href{mailto:#2}{#2}}}\frontmatter@RRAPformat}
{}{}
}%
\makeatother

\title{Stabilizing PDE--ML coupled systems}

\author[1]{Saad Qadeer\thanks{{\tt saad.qadeer@pnnl.gov}}} 
\author[2]{Panos Stinis} 
\author[3]{Hui Wan} 
\affil[1]{Advanced Computing, Mathematics and Data Division, Pacific Northwest National Laboratory, Richland, WA}
\affil[2]{Advanced Computing, Mathematics and Data Division, Pacific Northwest National Laboratory, Richland, WA}
\affil[ ]{Department of Applied Mathematics, University of Washington, Seattle, WA}
\affil[3]{Atmospheric Sciences and Global Change Division, Pacific Northwest National Laboratory, Richland, WA}

\begin{document}

\maketitle

\begin{abstract}
A long-standing obstacle in the use of machine-learnt surrogates with larger PDE systems is the onset of instabilities when solved numerically. Efforts towards ameliorating these have mostly concentrated on improving the accuracy of the surrogates or imbuing them with additional structure, and have garnered limited success. In this article, we study a prototype problem and draw insights that can help with more complex systems. In particular, we focus on a viscous Burgers'-ML system and, after identifying the cause of the instabilities, prescribe strategies to stabilize the coupled system. To improve the accuracy of the stabilized system, we next explore methods based on the Mori--Zwanzig formalism. 
\end{abstract}

\section{Introduction}\label{SecIntro}
Partial differential equations (PDEs) are an essential modeling tool in engineering and physical sciences. The numerical methods used for solving the more descriptive and sophisticated of these models comprise many  computationally expensive modules. Machine learning (ML) provides a way of replacing some of these modules by surrogates that are much more efficient at the time of inference. The resulting PDE--ML coupled systems, however, can be highly susceptible to instabilities \cite{brenowitz2018prognostic,brenowitz2020interpreting,ott2020fortran}. Efforts towards ameliorating these have mostly concentrated on improving the accuracy of the surrogates, imbuing them with additional structure, or introducing problem-specific stabilizers, and have garnered limited success \cite{rasp2018deep,wang2021stable,han2023ensemble,hu2024stable}. In this article, we study a prototype problem to understand the mathematical subtleties involved in PDE--ML coupling, and draw insights that can help with more complex systems.

\section{Stabilizing the coupled system}\label{SecStabilizing}

We study the viscous Burgers' equation with periodic boundary conditions
\eqn{
u_t  = \pz\left(-\frac{1}{2} u^2 \right) + \nu u_{zz},  \quad u(0,z) = u_0(z), \quad z\in [0,2\pi), \label{ViscBurg}
}
and consider the effect of replacing the diffusion term by a machine-learnt surrogate. We conduct online training by using the exact solutions for initial conditions of the form 
\eqn{
u_0(z) = \sum_{l= 0}^L a_l \cos(lz) + \sum_{l = 1}^L b_j \sin(lz), \qquad a_l,b_l \sim \frac{\text{unif}(-1,1)}{1+l^2}, \label{ICs}
}
to train the surrogate term. We first find that successful training necessarily requires the diffusion term to be written in divergence form, i.e.,
\eqn{
w_t = \pz\left(-\frac{1}{2}w^2 + \nu \mathcal{N}[w]\right), \label{MLSurrogate}
}
where $\mathcal{N}$ is a trainable architecture. Even then, the inability of the trained network to learn the high frequencies (the {\it spectral bias} phenomenon) means that naively solving the the resulting system \eqref{MLSurrogate} leads to instabilities for some initial conditions, or vastly incorrect solutions (relative errors of $O(1)$) for the rest (even within the temporal training interval). 

These observations can be explained by noting that the spectral bias limits the ability of \eqref{MLSurrogate} to correctly treat the high frequencies, in effect leading to a lopsided under-resolved model, since the high frequency modes keep picking up energy due to the convective term. To ameliorate this, we include a low-pass filter in this term, to wit, 
\eqn{
\phi_t = \pz\left(-\frac{1}{2}\mathcal{P}_M\left[\phi^2\right] + \nu \mathcal{P}_M\left[\mathcal{N}[\phi]\right]\right), \label{MLSurrLP}
}
where
\eqn{
\mathcal{P}_M\left[g\right](z) = \sum_{|j| \leq M} \widehat{g}_je^{ijz}  \text{ for } g(z) = \sum_{j \in \mathbb{Z}} \widehat{g}_je^{ijz}. \label{LPFilter}
}

\begin{figure}[tbph]
\centering
\subfigure[$u_0(z) = \sin(z)$]
{\includegraphics[width=0.47\textwidth]{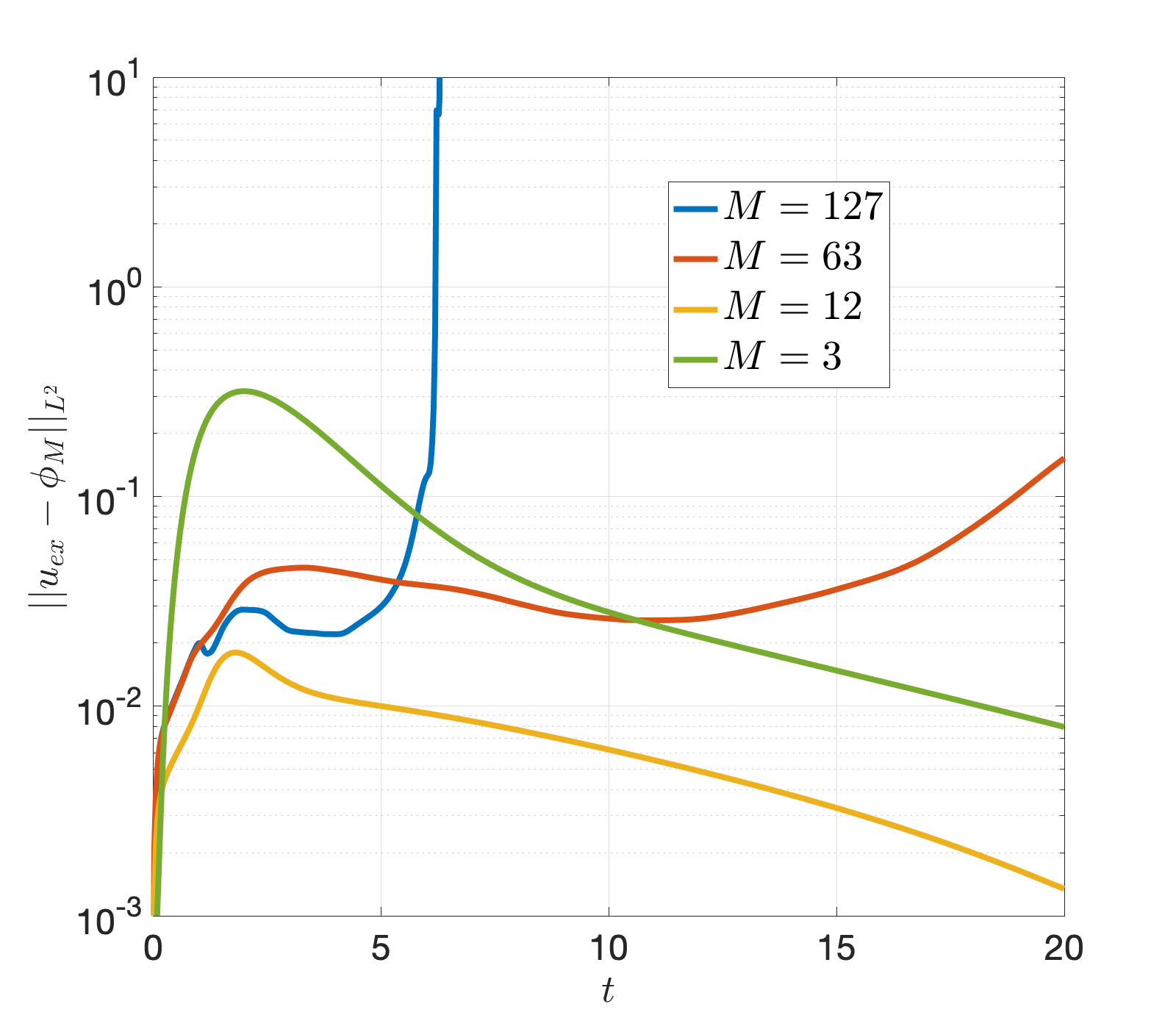}\label{SpecBiasSin}
}
\subfigure[$u_0(z) = e^{\sin(z)}$]
{\includegraphics[width=0.47\textwidth]{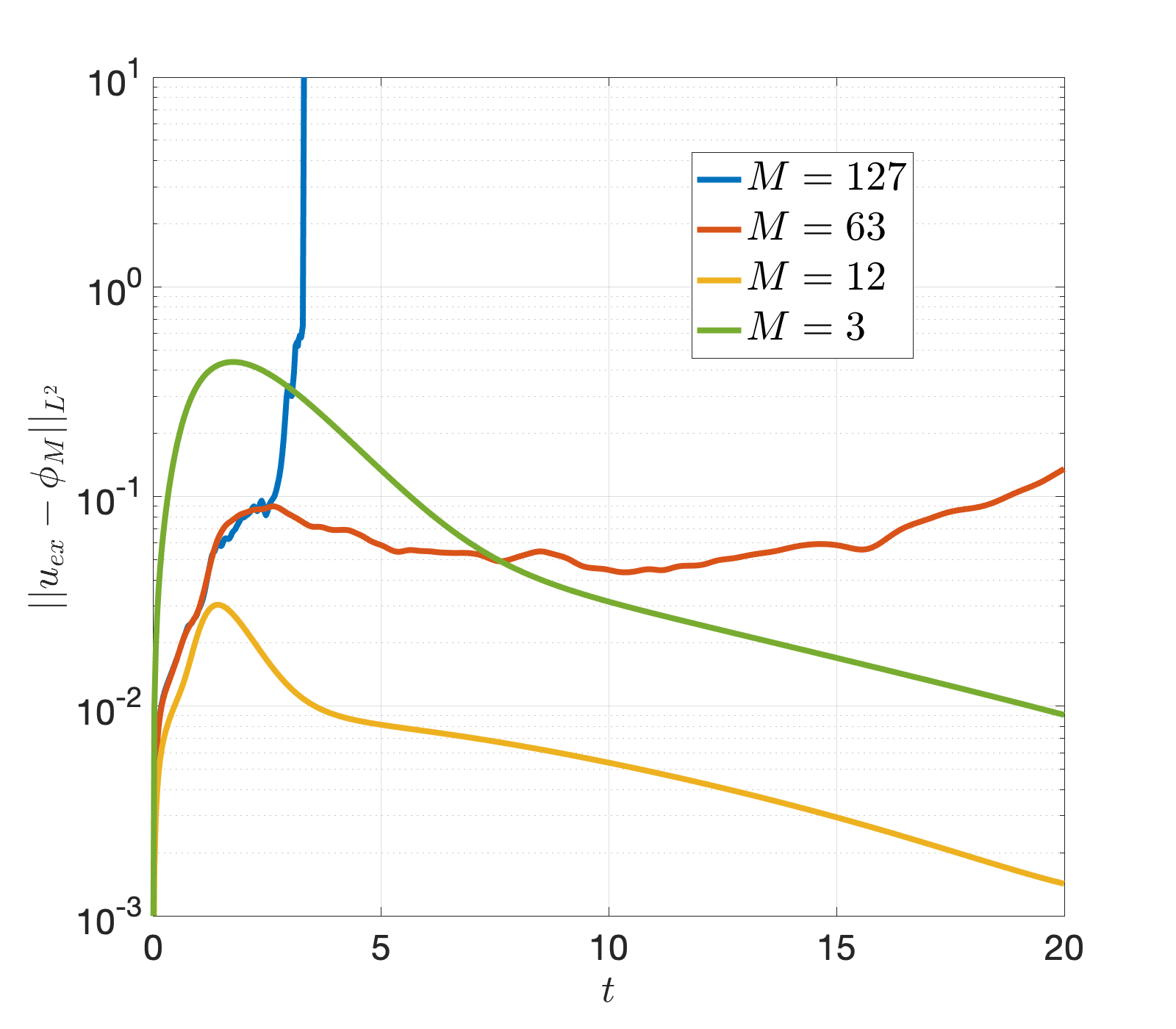}\label{SpecBiasExpSin}
}
\caption{The $L^{2}$ error evolution with time for different values of $M$, the cut-off for the low-pass filter.}\label{SpecBias}
\end{figure}

In Figure \ref{SpecBias}, we show the results for various values of $M$. The PDEs are solved on a uniform grid of size $2^8$; the surrogate network has widths $2^8$, $2^9$, and $2^8$ (so only two layers of parameters) and is trained on solution snapshots over $[t_i,t_i+10\Delta t] \subset [0,4]$ (50 such initial conditions, 250 sets of such snapshots each) with $L = 24$. We see that $M = 127$ (no band-limiting) gives either $O(1)$ errors or instabilities, whereas $M = 12$ improves both stability and accuracy. (We also use the exponentially decaying spectral filter to avoid aliasing as usual.) Lowering this further to $M = 3$ leads to stable but somewhat inaccurate results. 

To sum up, we find that \eqref{MLSurrLP} yields stable solutions for all initial conditions that are more accurate than those for \eqref{MLSurrogate}. Since it is more closely aligned with a Galerkin (i.e., Markovian) treatment of \eqref{ViscBurg} using a small number of modes, we expect that introducing memory-based corrections can help improve the accuracy of the stabilized system.

\section{Introducing memory terms}\label{SecMemoryLap}
Given a system of ODEs
\eqn{
\phi'_k(t) = R_k({\bs \phi}), \quad \phi_k(0) = x_k, \qquad k = 1,\hdots,N, \label{ODEsys}
}
we wish to construct a reduced-order model (ROM) for the ``resolved'' variables $\widehat{\bs{\phi}} = \left(\phi_l\right)_{1 \leq l \leq N_\text{rom}}$ for some $N_\text{rom} < N$ that reproduces the dynamics of the full system. Let $\mathcal{L} = \sum_{j = 1}^N R_j(x) \frac{\partial}{\partial x_j} $ be the Liouville operator corresponding to \eqref{ODEsys}. Define the projection operator $P$ by $[Pg](\widehat{x}) = g(\widehat{x},0)$ for any function $g(\widehat{x},\widetilde{x})$ of both $\widehat{x}$, the resolved, and $\widetilde{x}$, the unresolved initial conditions. Further set $Q = I - P$ and
\eqn{
\widecheck{g}(\widehat{x}) = [Pg](\widehat{x}), \qquad \overline{g}(x) = [Qg](x) = g(x) - [Pg](\widehat{x}). \nonumber
}

The Mori--Zwanzig (MZ) formalism then states that \cite{chorin2007problem}
\eqn{
\phi_k'(t) = e^{t\mathcal{L}}\mathcal{L}x_k = e^{t\mathcal{L}}P\mathcal{L}x_k + \int_0^t e^{(t-s)\mathcal{L}}P\mathcal{L}e^{sQ\mathcal{L}}Q\mathcal{L}x_k \ ds +  e^{tQ\mathcal{L}}Q\mathcal{L}x_k. \label{MZId}
}

In particular, defining $\bs{\eta}(t;x) = \begin{pmatrix}\eta_1 & \hdots & \eta_N\end{pmatrix}$ as the solution to the orthogonal dynamics equations 
\eqn{
\eta_l'(t) = \overline{R}_l(\bs{\eta}), \qquad \bs{\eta}(0) = x, \qquad 1 \leq l \leq N,\label{ODEqns}
} 
allows us to write, for $1 \leq k \leq N_\text{rom}$, 
\eqn{
\phi'_k(t ; x) = \widecheck{R}_k\left(\widehat{\phi}(t ; x)\right) + \int_0^t e^{(t-s)\mathcal{L}}P\mathcal{L}e^{sQ\mathcal{L}}Q\mathcal{L}x_k \ ds +  \overline{R}_k\left(\bs{\eta(t;x)}\right), \label{MZform1}
}

The integral in \eqref{MZform1} is termed the memory and the third term is called the noise; observe that $e^{sQ\mathcal{L}}Q\mathcal{L}x_k = \overline{R}_k\left(\bs{\eta(s;x)}\right)$ so both the memory integrand and the noise rely on the solution to the orthogonal dynamics. Because of this, the system \eqref{MZform1} is not closed. 

However, we can use \eqref{MZform1} to construct ROMs under carefully justified assumptions and approximations. First, by limiting ourselves to initial conditions of the form $x = (\widehat{x},0)$ only, we find that $\bs{\eta}(t;\widehat{x},0) = (\widehat{x},0)$, and hence
\eqn{
\overline{R}_k\left(\bs{\eta(t;\widehat{x},0)}\right) = \overline{R}_k\left(\widehat{x},0\right) = 0, \nonumber
}
so the noise term vanishes. Observe that limiting ourselves to these initial conditions is equivalent to applying $P$ throughout \eqref{MZId}. We therefore obtain
\eqn{
\phi'_k(t ; \widehat{x},0) = \widecheck{R}_k\left(\widehat{\phi}(t ; \widehat{x},0)\right) + \int_0^t Pe^{(t-s)\mathcal{L}}P\mathcal{L}e^{sQ\mathcal{L}}Q\mathcal{L}x_k \ ds. \label{MZform2}
}

We note that this assumption does not help simplify the solution of the orthogonal dynamics in the memory integrand as the application of the Liouville operator to it entails considering the rate of change with respect to the initial conditions.

Following \cite{chorin2002optimal}, we develop an approach that allows us to rewrite \eqref{MZform2}  in terms of memory kernels by relying on an ensemble of solutions generated according to a well-defined rule. Define the weight function
\eqn{
w(x) = \Pi_{i = 1}^{N_\text{rom}} \frac{e^{-\frac{(x_i-\mu_i)^2}{2\sigma_i^2}}}{\sigma_i\sqrt{2\pi}} \Pi_{j = N_\text{rom}+1}^N \delta_0(x_j), \label{WeightFcn}
}
and denote
\eqn{
\ip{g_1,g_2} = \int_{\mathbb{R}^N} g_1(x)g_2(x) w(x) \ dx. \label{IPDefn}
}

The projection $P$ can be characterized as the conditional expectation
\eqn{
\left[Pg\right](\widehat{x}) = \frac{\int_{\mathbb{R}^{N-N_\text{rom}}} g(\widehat{x},\widetilde{x}) w(x) \ d\widetilde{x}}{\int_{\mathbb{R}^{N-N_\text{rom}}} w(x) \ d\widetilde{x}} = g(\widehat{x},0), \nonumber
}
as it simply sets all the unresolved variables to zero. It can also be viewed as an orthogonal projection with respect to $\ip{\cdot,\cdot}$ on the subspace of all functions of $\widehat{x}$ since, for any function $h$ in this subspace, 
\eqn{
\ip{g-Pg , h} = \int_{\mathbb{R}^N} \left[g(x) - g(\widehat{x},0)\right]h(\widehat{x}) w(x) \ dx = \int_{\mathbb{R}^{N_\text{rom}}} \left[g(\widehat{x},0) - g(\widehat{x},0)\right]h(\widehat{x}) w(\widehat{x},0) \ d\widehat{x} = 0. \nonumber
}

Let $\{h_j(\widehat{x})\}_{j \geq 1}$ be an orthonormal basis for the range of $P$. We can then write
\eqn{
[Pg](\widehat{x}) = \sum_{j \geq 1} \ip{h_j , g} h_j(\widehat{x}). \label{HermProj}
}

Next, we use this ideas to simplify the integrand in \eqref{MZform2}. Applying Dyson's identity to $Q\mathcal{L}x_k$ leads to
\eqn{
e^{tQ\mathcal{L}}Q\mathcal{L}x_k = e^{t\mathcal{L}}Q\mathcal{L}x_k - \int_0^t e^{ (t-s)\mathcal{L}} P\mathcal{L} e^{sQ\mathcal{L}}Q\mathcal{L}x_k \ ds. \label{DysonId} 
}

Defining $A_k(t;x) = e^{tQ\mathcal{L}} Q\mathcal{L}x_k$ and noting then that $A_k(0;x) = Q\mathcal{L}x_k$, we obtain
\eqn{
A_k(t ; x) = e^{t\mathcal{L}}A_k(0 ; x) - \int_0^t e^{ (t-s)\mathcal{L}} P\left(\mathcal{L} A_k(s ; x) \right) \ ds. \nonumber 
}

Applying $P\mathcal{L}$ throughout gives
\eqn{
P\left(\mathcal{L}A_k(t ; x)\right) = P\left(\mathcal{L}e^{t\mathcal{L}}A_k(0 ; x)\right) - P \int_0^t \mathcal{L}e^{ (t-s)\mathcal{L}} P\left(\mathcal{L} A_k(s ; x) \right) \ ds. \label{PLeta}
}

Expand
\eqn{
\left[P\mathcal{L}A_k\right](t;\widehat{x}) = \sum_{j \geq 1} \ip{h_j , \mathcal{L}A_k(t)} h_j(\widehat{x}) =  \sum_{j \geq 1}K_{jk}(t)  h_j(\widehat{x}) \label{MemKerExp}
}
where $K_{jk}(t) := \ip{h_j , \mathcal{L}A_k(t)}$ are termed memory kernels. Plugging this expansion in \eqref{PLeta} and taking inner products with $h_l$ throughout yields
\eqn{
\ip{h_l , \mathcal{L} A_k(t)} = \ip{h_l , \mathcal{L}e^{t\mathcal{L}} A_k(0)} - \sum_{j \geq 1} \int_0^t \ip{h_l , \mathcal{L}e^{(t-s)\mathcal{L}}h_j} K_{jk}(s) \ ds, \nonumber
}
that is,
\eqn{
K_{lk}(t) = f_{lk}(t) - \sum_{j \geq 1} \int_0^t g_{lj}(t-s) K_{jk}(s) \ ds, \label{KInteqn}
}
where
\eqn{
f_{lk}(t) := \ip{h_l , e^{t\mathcal{L}} \mathcal{L} A_k(0)}, \quad g_{lj}(t) :=  \ip{h_l , e^{t\mathcal{L}}\mathcal{L}h_j}. \label{fgdefs}
}

Observe that we have used the commutativity of $e^{t\mathcal{L}}$ and $\mathcal{L}$ in \eqref{fgdefs} to simplify the resulting calculations; we also note that in the first expression, since $A_k(0;x) = \overline{R}_k(x)$, the operator $\mathcal{L}$ acts on functions of $x$, and not just $\widehat{x}$; in the second, since $h_j$ is a function of $\widehat{x}$, the operator $\mathcal{L}$ acts on only these but also yields a function of $x$. The corresponding evolution operators then act on all the variables so we need to track the trajectories of all the variables in the full-order system; the projections implicit in the inner products in \eqref{fgdefs} then have the effect that we only need to consider initial conditions with the unresolved variables set to zero. Discretizing the convolution integral in \eqref{KInteqn} by, e.g., the trapezoidal rule then yields a linear system that we can solve at each timestep to compute $K_{lk}(t_n)$ (observe that $\left(g_{lj}(t)\right)$ is a square matrix for each $t$). 

The evaluation of the inner products can be carried out either by Monte--Carlo integration, i.e., by sampling many points according to the weight function $w(x)$ and averaging, or by employing an accurate quadrature rule. The particular form of $w(x)$ makes the latter approach more appealing. Specifically, we note that
\eqn{
\int_{-\infty}^{\infty} g(x) \frac{e^{-\frac{(x - \mu)^2}{2\sigma^2}}}{\sigma \sqrt{2\pi}} \ dx = \frac{1}{\sqrt{2\pi}} \int_{-\infty}^{\infty} g\left(\sigma z + \mu \right) e^{-\frac{z^2}{2}} \ dz \approx \sum_{i = 1}^{N_q} g(\sigma z_i + \mu) \omega_i, \label{GHQuad}
}
where $\{(z_i,\omega_i)\}_{1 \leq i \leq N_q}$ is the Gauss--Hermite quadrature rule, exact up to polynomials of degree $(2N_q - 1)$. Taking tensor products of this one-dimensional rule for all the resolved variables yields the high-dimensional analog. If the same sized grid is used for each resolved variable, we end up requiring $(N_q)^{N_\text{rom}}$ full-order solutions. With that said, the solution trajectories can be calculated in parallel and used to compute the $\left\{f_{lk}(t)\right\}$ and $\left\{g_{lj}(t)\right\}$.

In a similar spirit, the orthonormal basis $\{h_j(\widehat{x})\}_{j \geq 1}$ can be taken to be tensor products of one-dimensional Hermite polynomials $\{H_j(z)\}$, corresponding to weight function $e^{-z^2/2}$, with the arguments modified suitably to account for the $\mu_j$ and $\sigma_j$. In practice, we can only employ a finite basis; if we allow polynomials up to highest degree $d$, we need a total of $J = {N_\text{rom}+d \choose d}$ polynomials. The resulting reduced-order model is
\eqn{
\phi'_k(t ; \widehat{x},0) = \widecheck{R}_k\left(\widehat{\phi}(t ; \widehat{x},0)\right) + \sum_{j = 1 }^J \int_0^t h_j\left(\widehat{\phi}(s ; \widehat{x},0)\right) K_{jk}(t-s) \ ds. \label{MemKernels2}
}

\begin{figure}[tbph]
\centering
\subfigure[$x_1 = 1$]
{\includegraphics[width=0.24\textwidth]{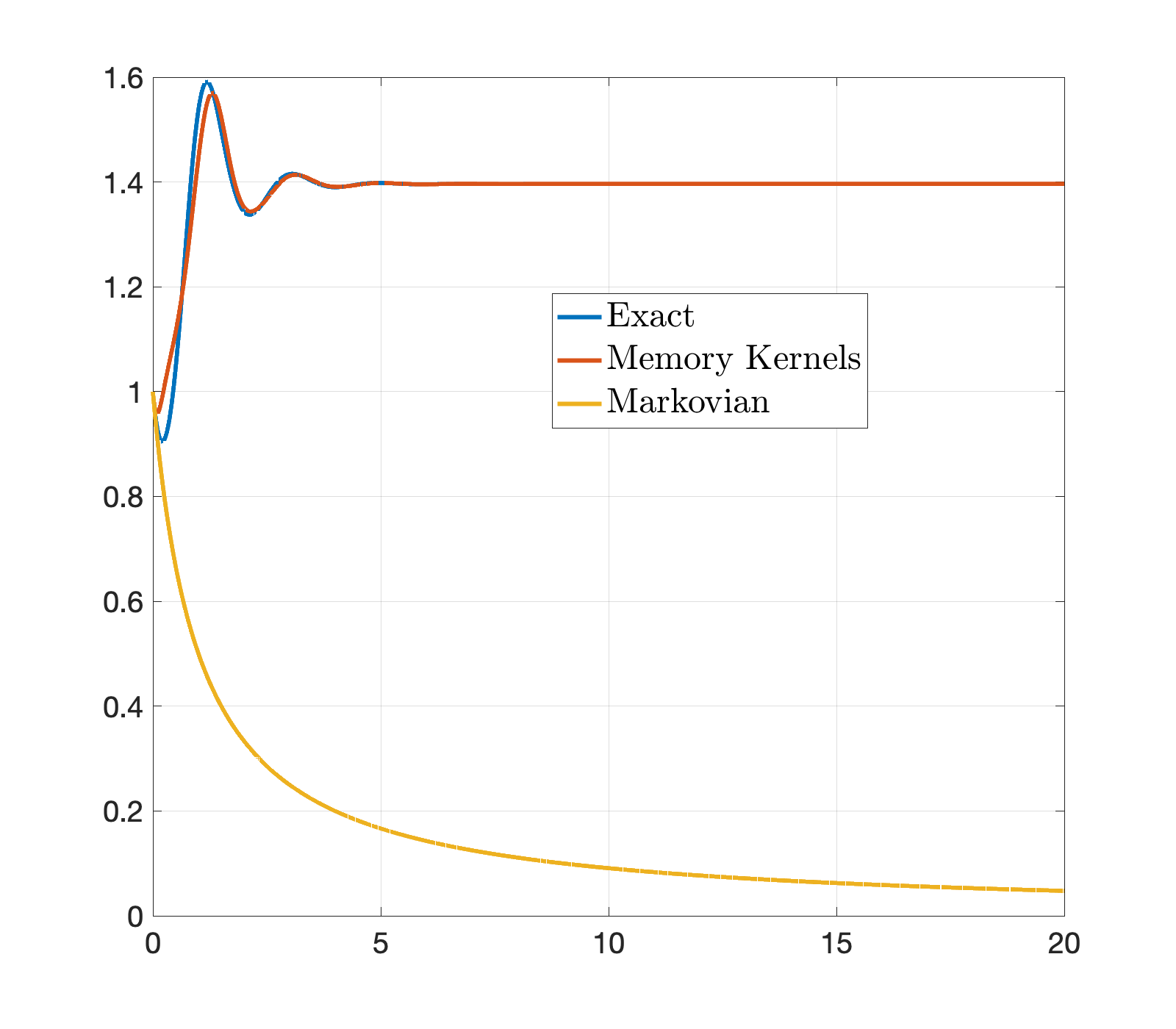}\label{HeQu1_1}
}
\subfigure[$x_1 = 2$]
{\includegraphics[width=0.24\textwidth]{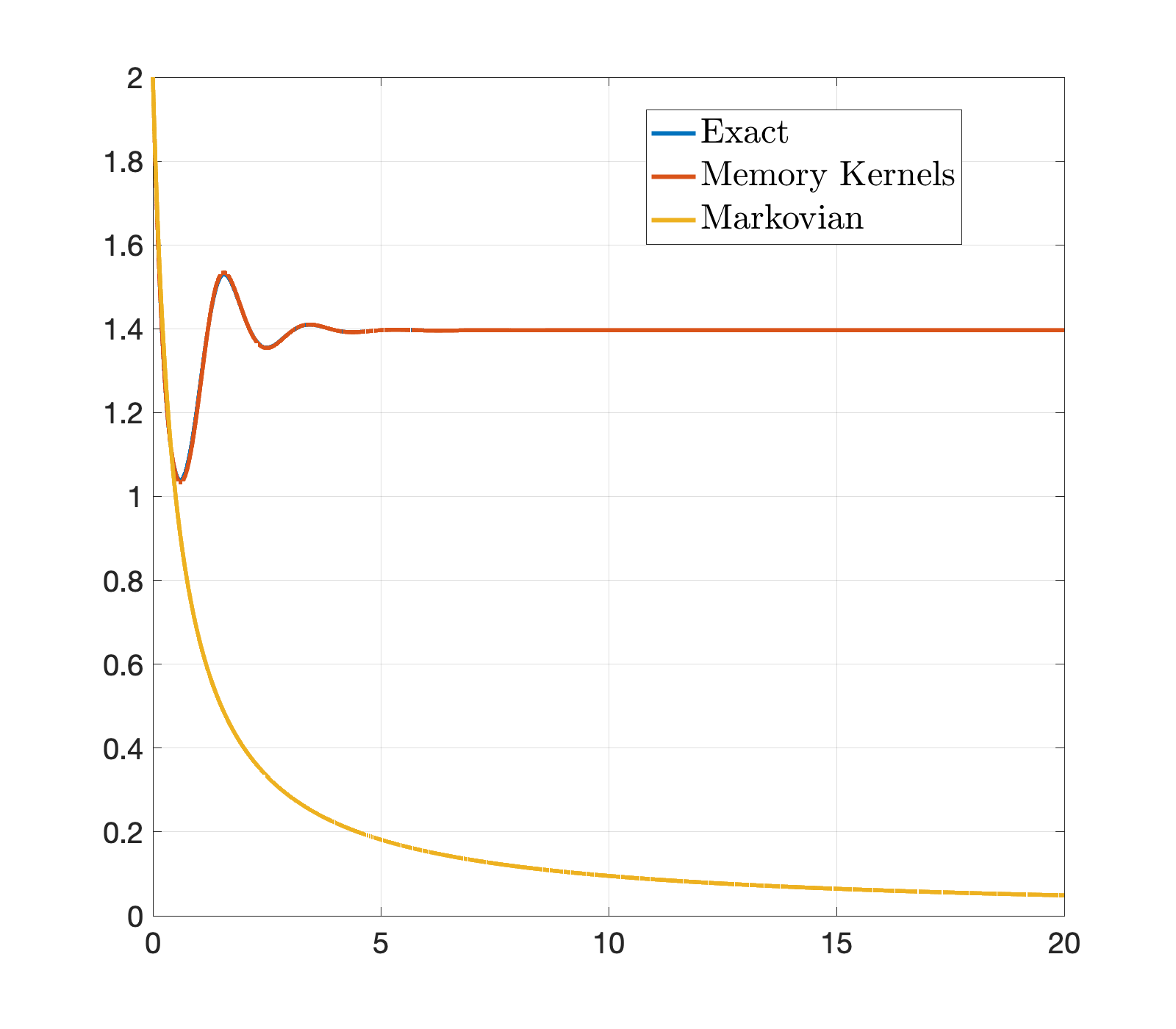}\label{HeQu1_2}
}
\subfigure[$x_1 = 3$]
{\includegraphics[width=0.24\textwidth]{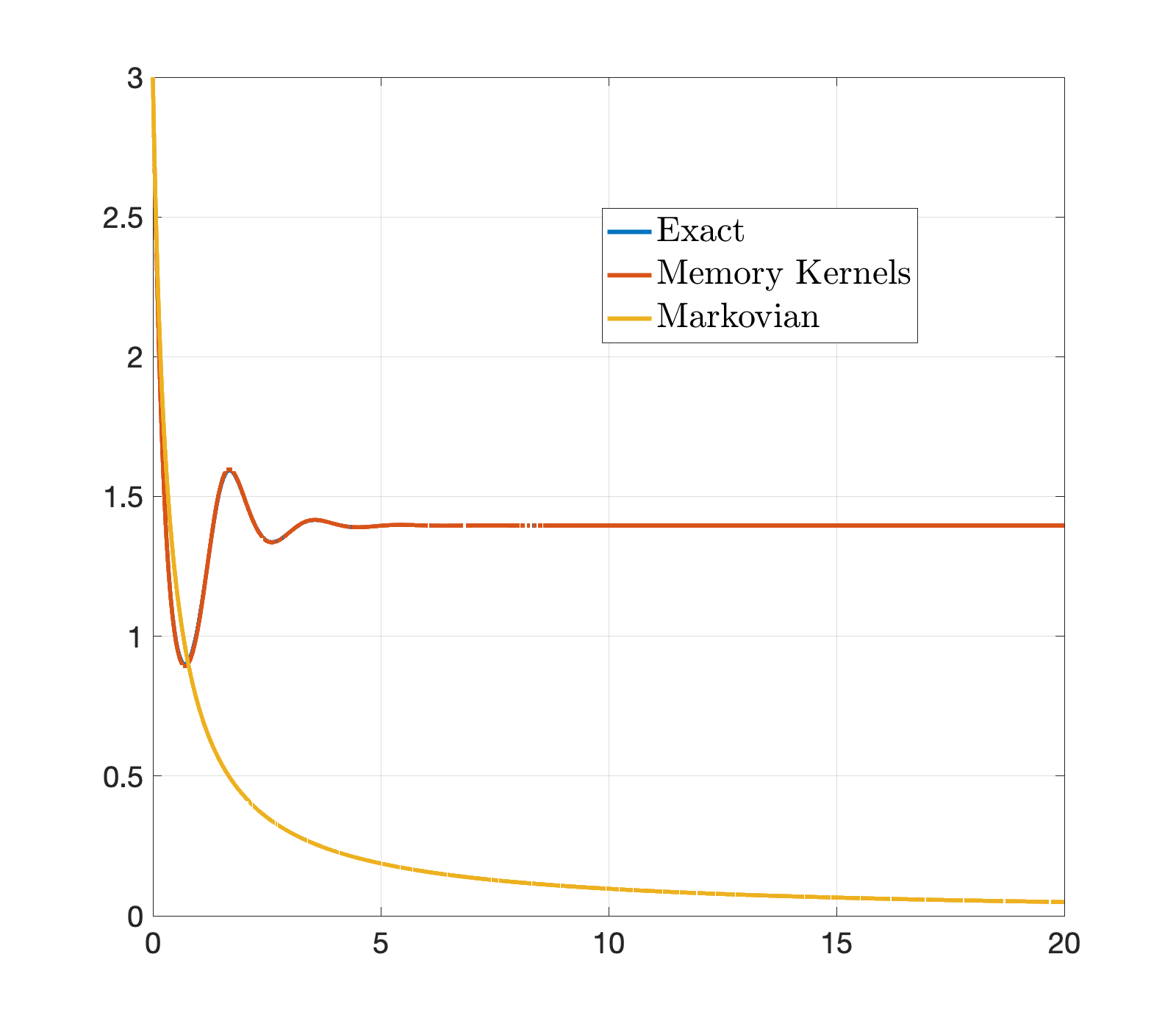}\label{HeQu1_3}
}
\subfigure[$x_1 = 7$]
{\includegraphics[width=0.24\textwidth]{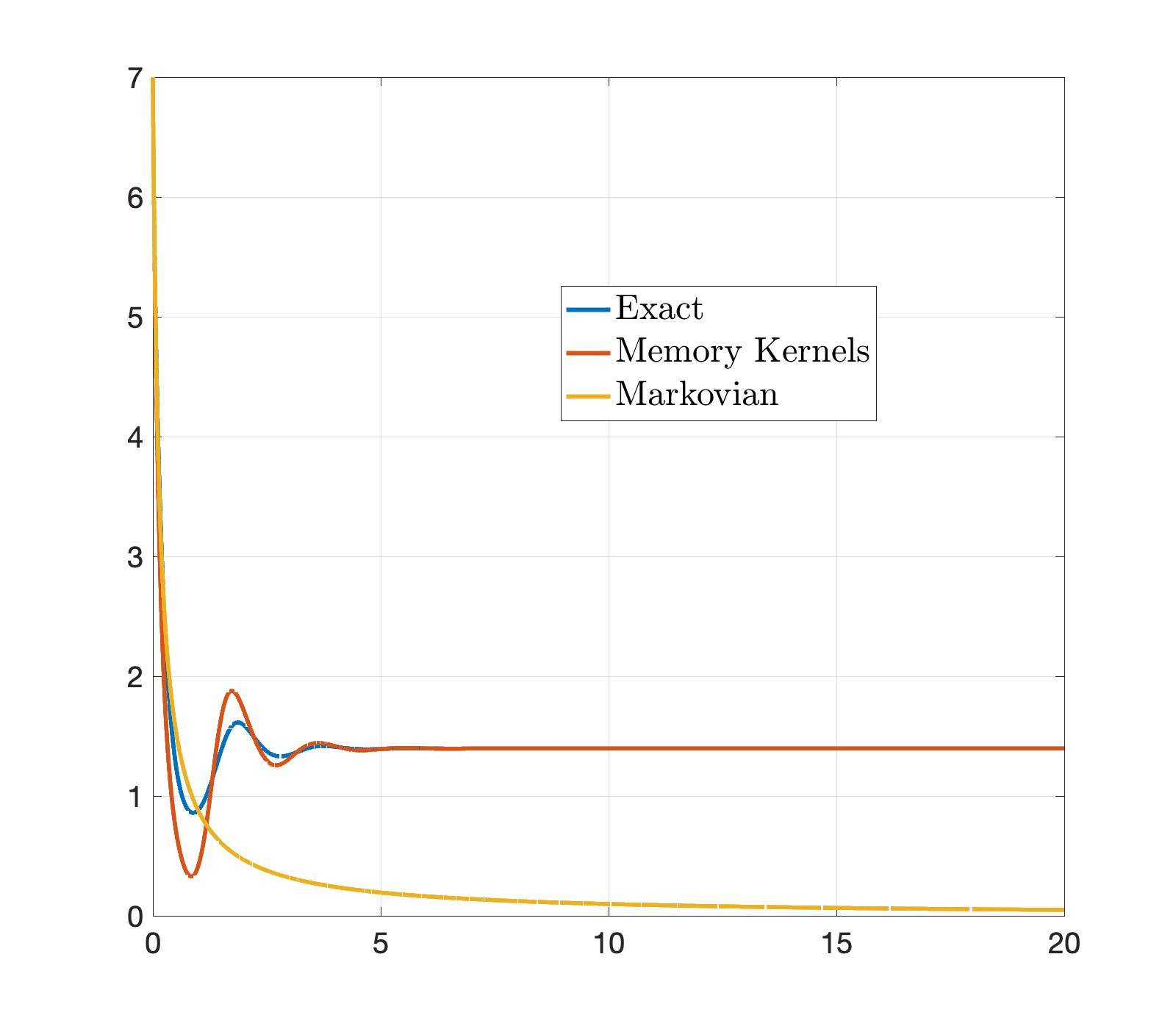}\label{HeQu1_7}
}
\caption{Comparison of the exact solution to \eqref{NLexample} and those obtained by expanding the memory kernels up to linear polynomials. We employ $N_q = 20$ Hermite quadrature nodes translated to $[1,4]$.}\label{HeQu1}
\end{figure}

\begin{figure}[tbph]
\centering
\subfigure[$x_1 = 1$]
{\includegraphics[width=0.24\textwidth]{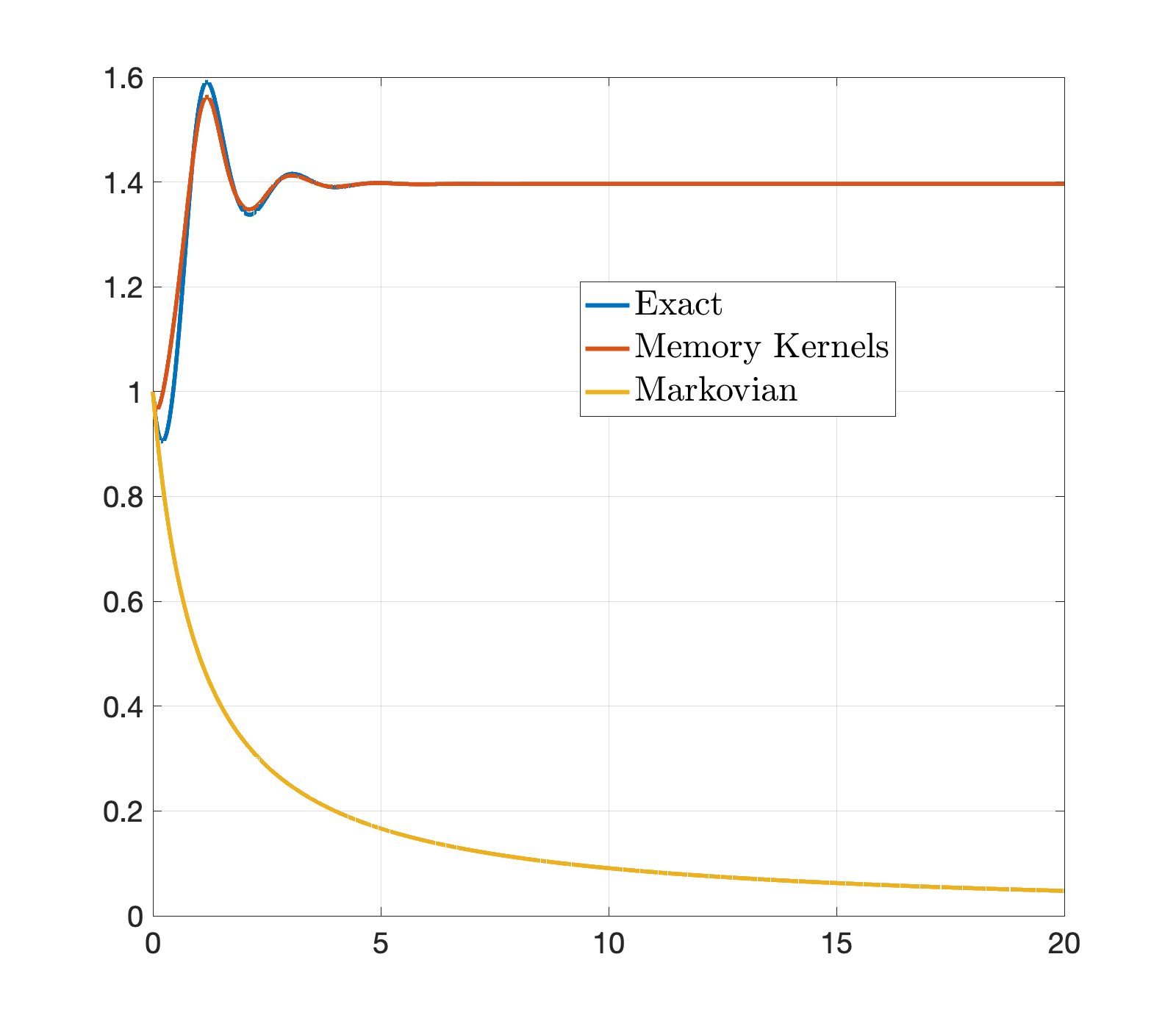}\label{HeQu2_1}
}
\subfigure[$x_1 = 2$]
{\includegraphics[width=0.24\textwidth]{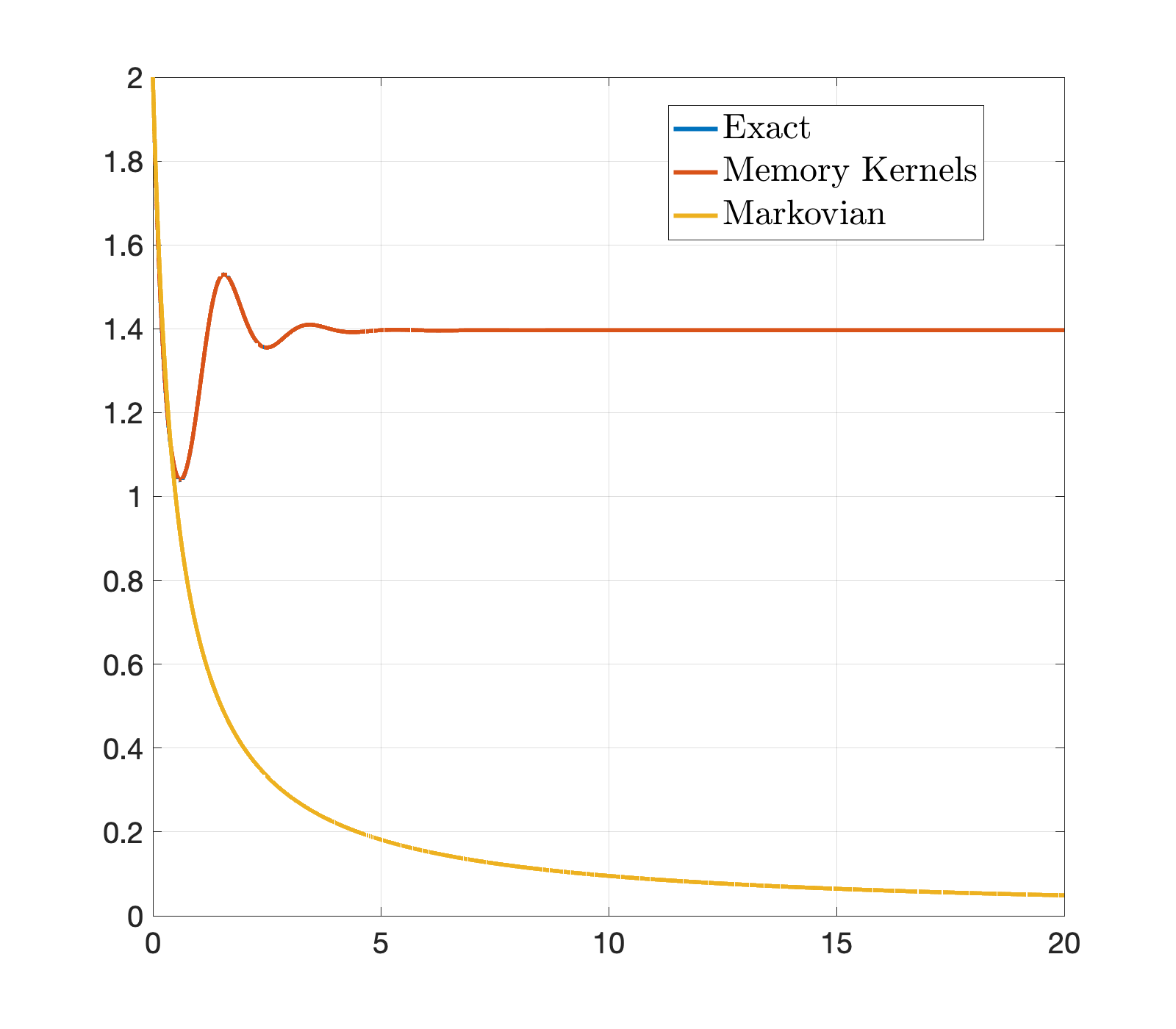}\label{HeQu2_2}
}
\subfigure[$x_1 = 3$]
{\includegraphics[width=0.24\textwidth]{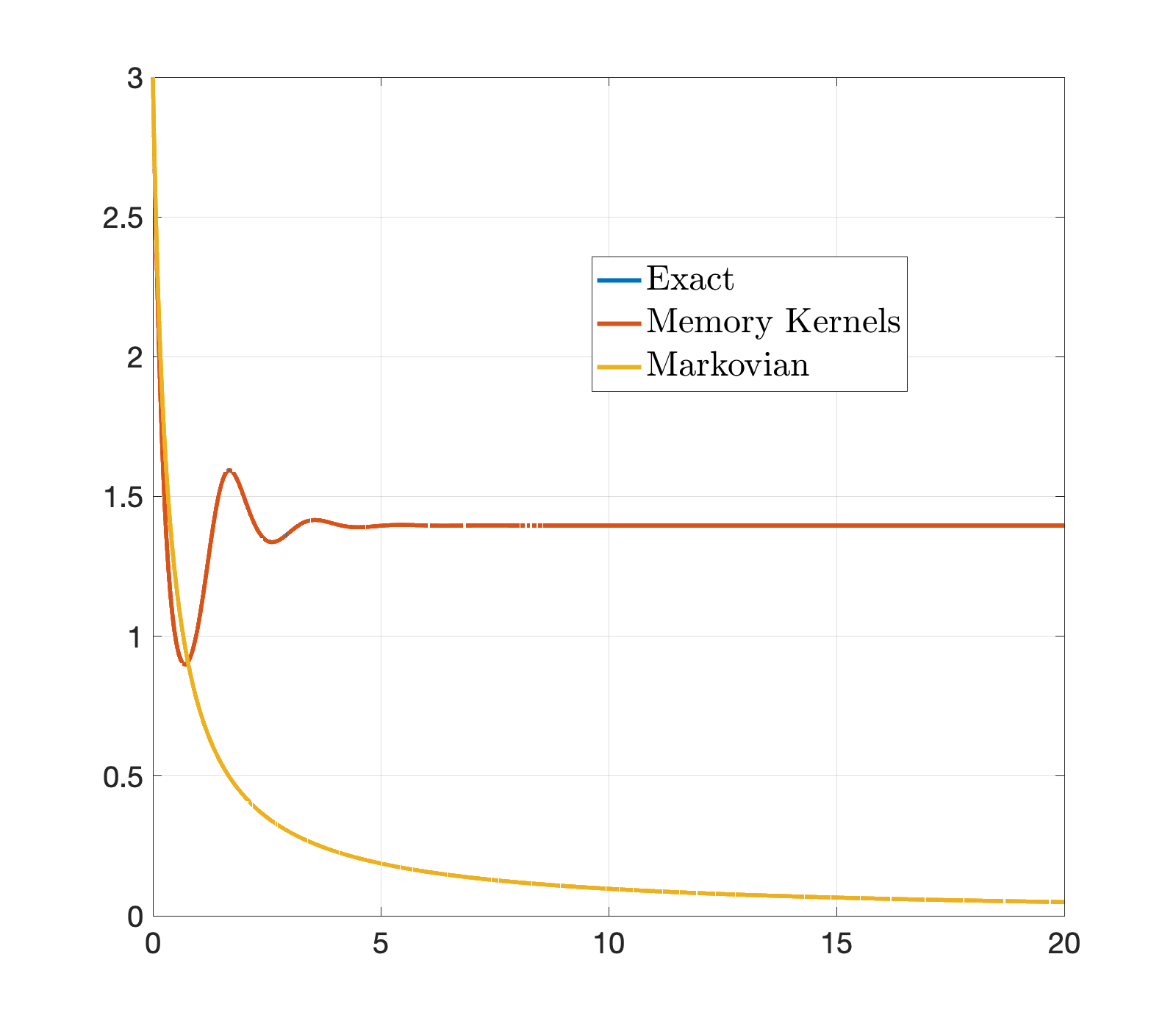}\label{HeQu2_3}
}
\subfigure[$x_1 = 7$]
{\includegraphics[width=0.24\textwidth]{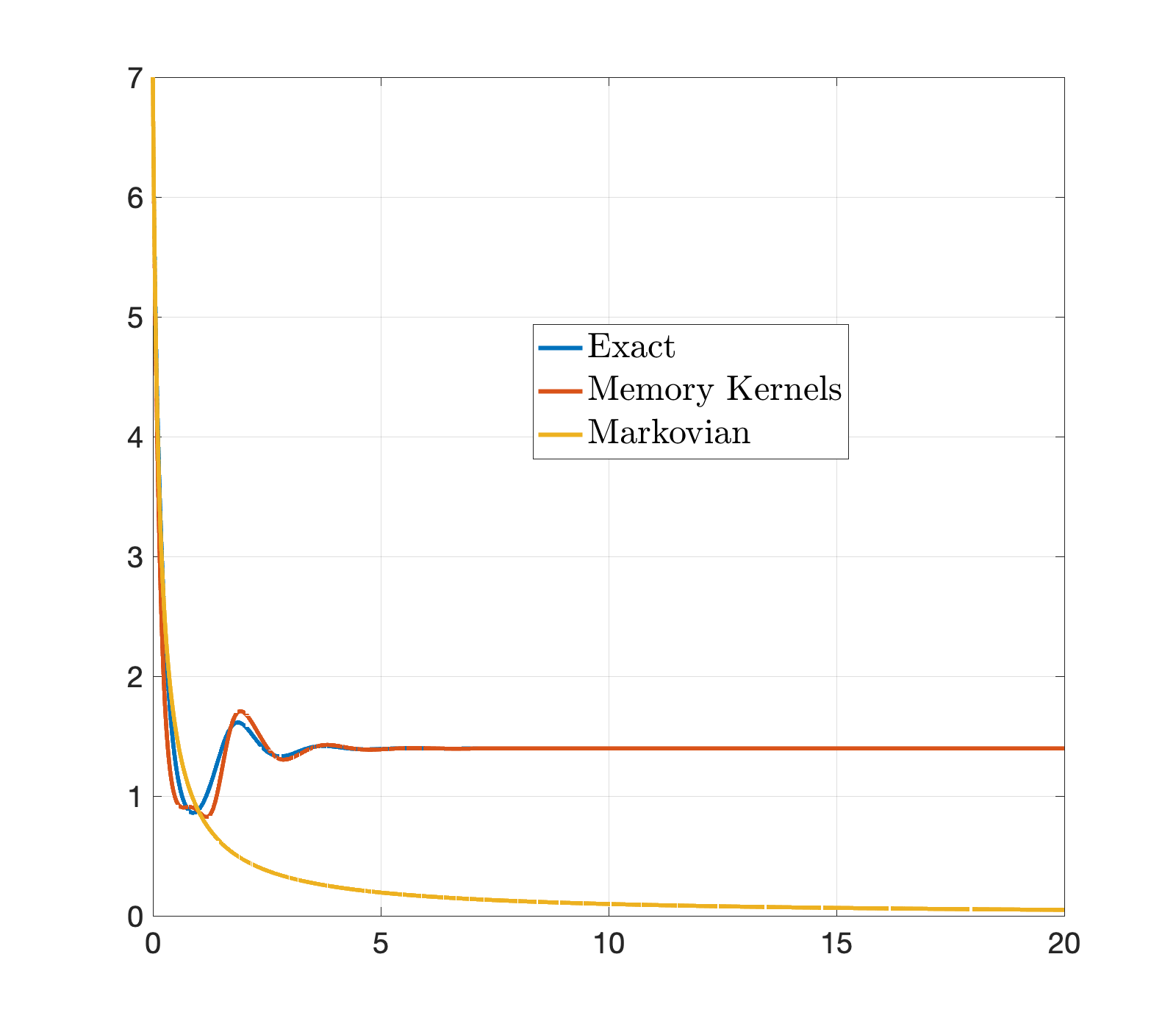}\label{HeQu2_7}
}
\caption{Comparison of the exact solution to \eqref{NLexample} and those obtained by expanding the memory kernels up to quadratic polynomials. We employ $N_q = 30$ Hermite quadrature nodes spaced over $[1,4]$.}\label{HeQu2}
\end{figure}

\begin{figure}[tbph]
\centering
\subfigure[$x_1 = 1$]
{\includegraphics[width=0.24\textwidth]{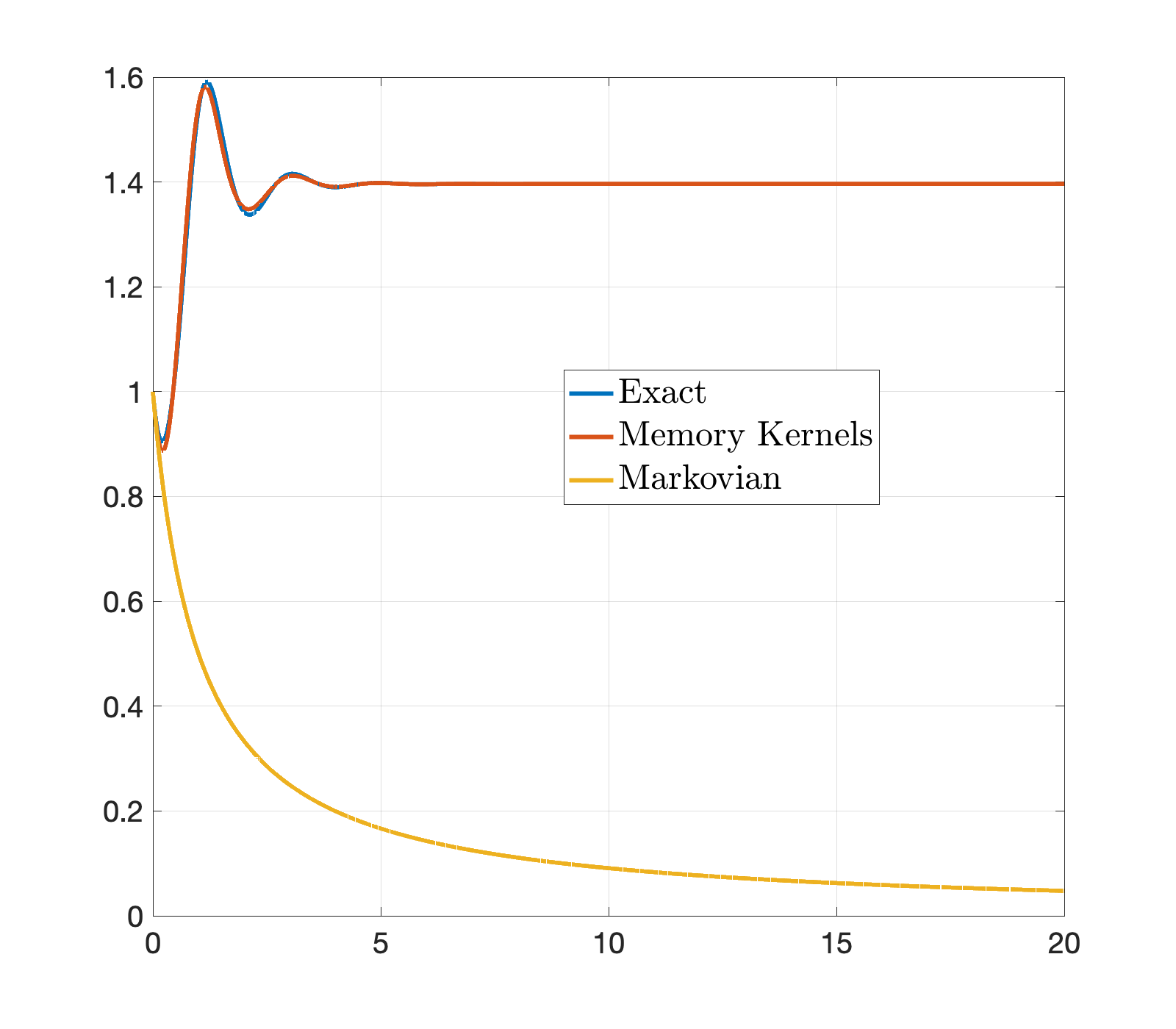}\label{HeQu3_1}
}
\subfigure[$x_1 = 2$]
{\includegraphics[width=0.24\textwidth]{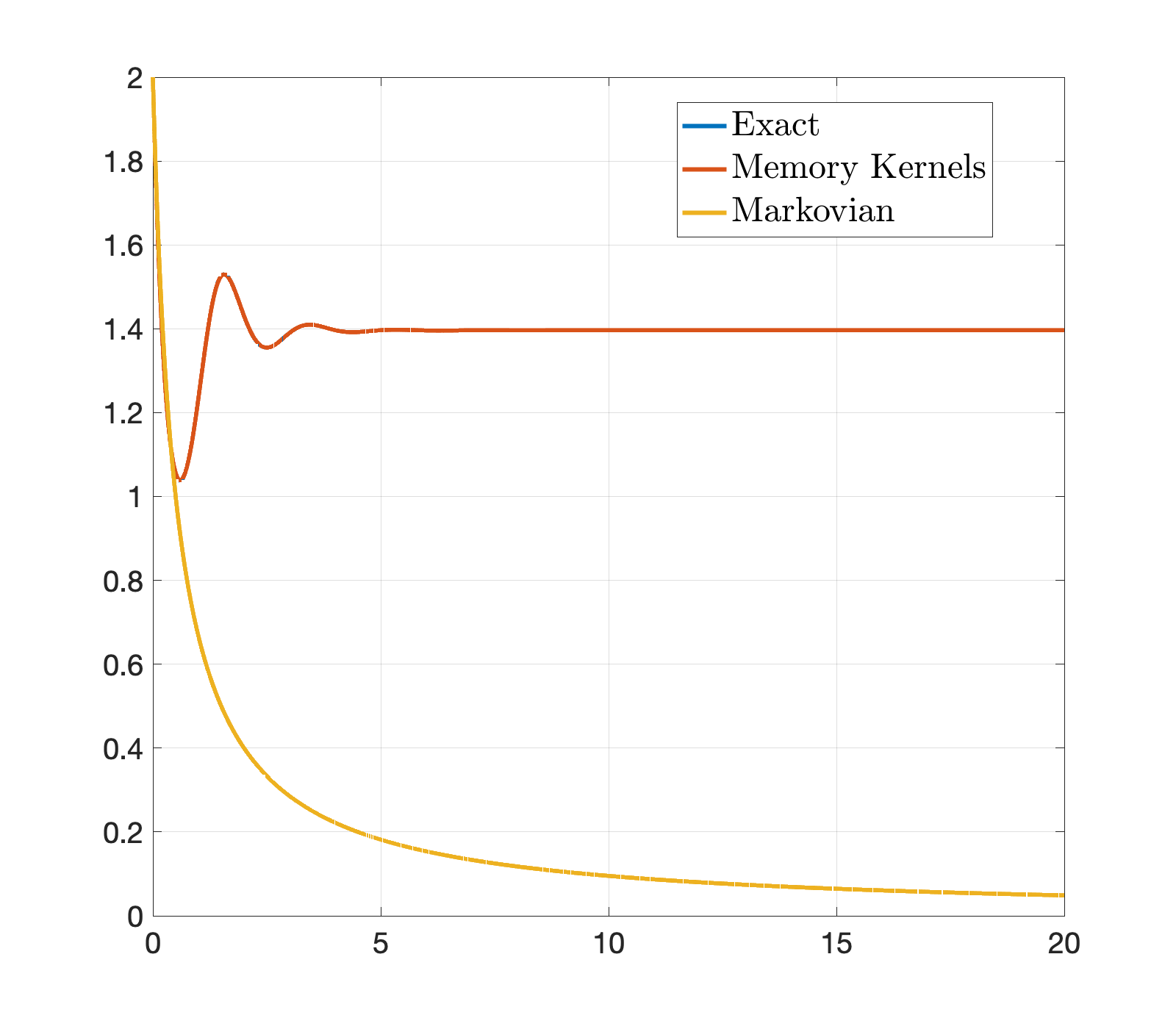}\label{HeQu3_2}
}
\subfigure[$x_1 = 3$]
{\includegraphics[width=0.24\textwidth]{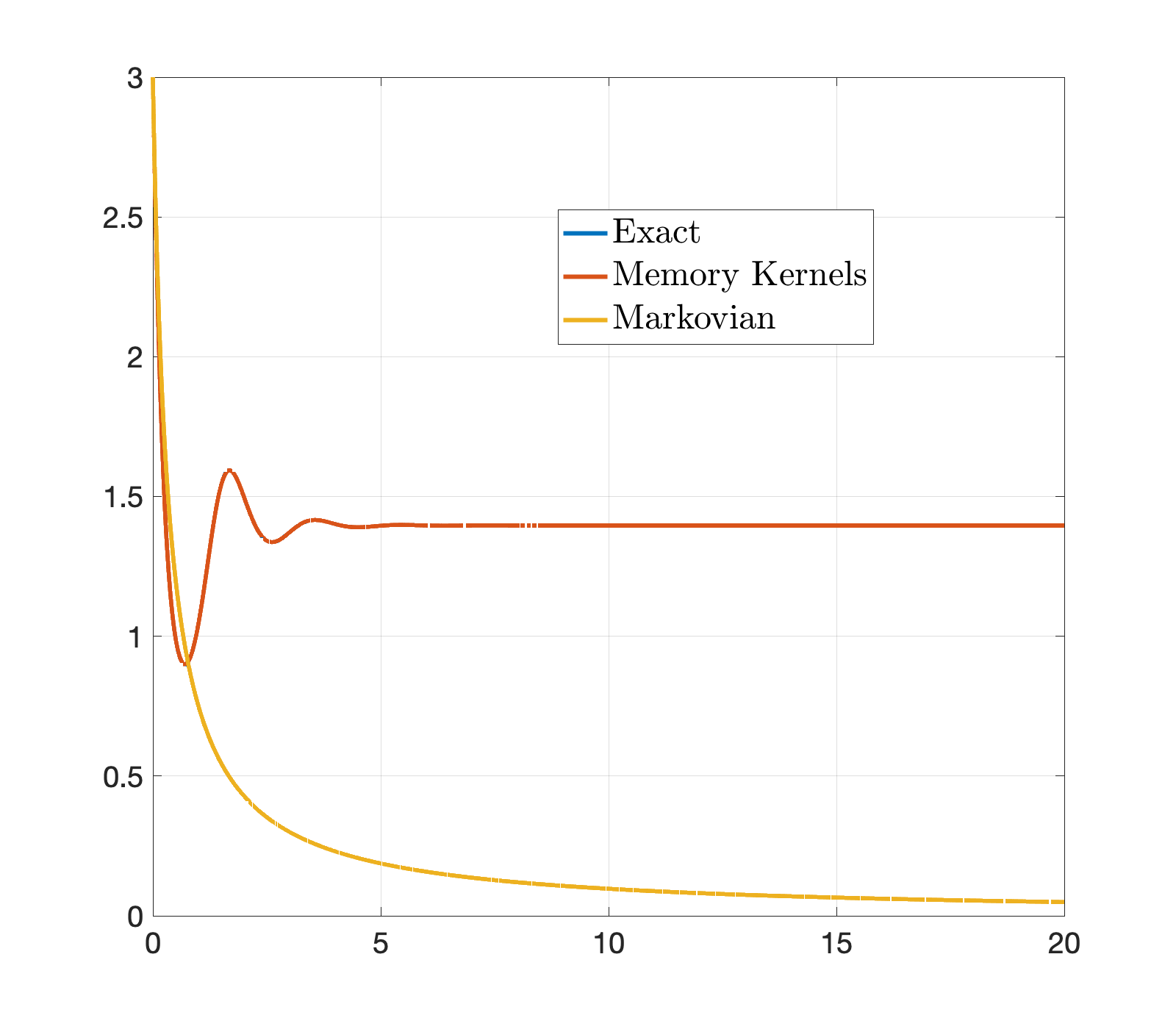}\label{HeQu3_3}
}
\subfigure[$x_1 = 7$]
{\includegraphics[width=0.24\textwidth]{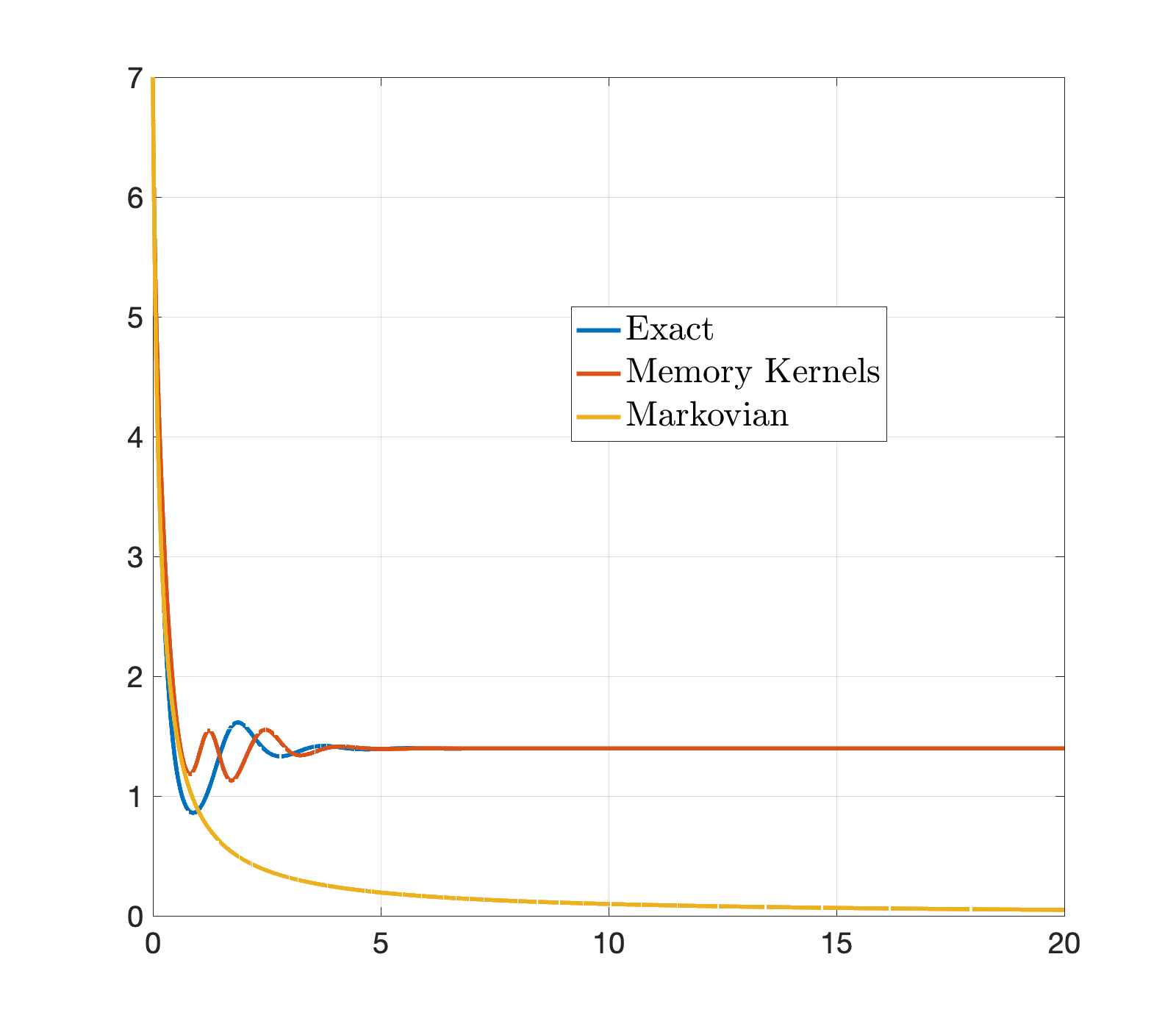}\label{HeQu3_7}
}
\caption{Comparison of the exact solution to \eqref{NLexample} and those obtained by expanding the memory kernels up to cubic polynomials. We employ $N_q = 40$ Hermite quadrature nodes translated to $[1,4]$.}\label{HeQu3}
\end{figure}

\begin{figure}[tbph]
\centering
\subfigure[$d = 1$]
{\includegraphics[width=0.31\textwidth]{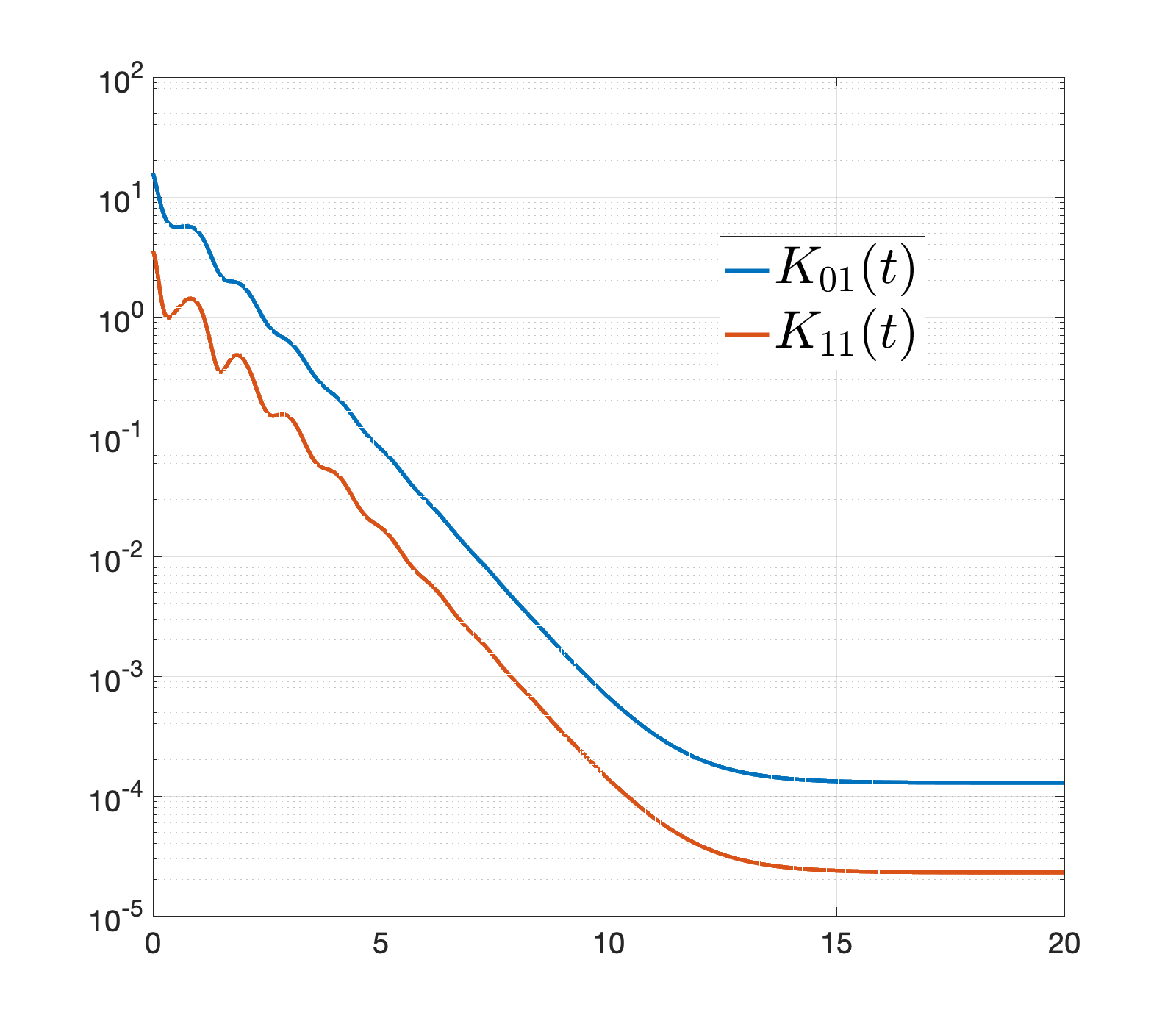}\label{HeQuK1}
}
\subfigure[$d = 2$]
{\includegraphics[width=0.31\textwidth]{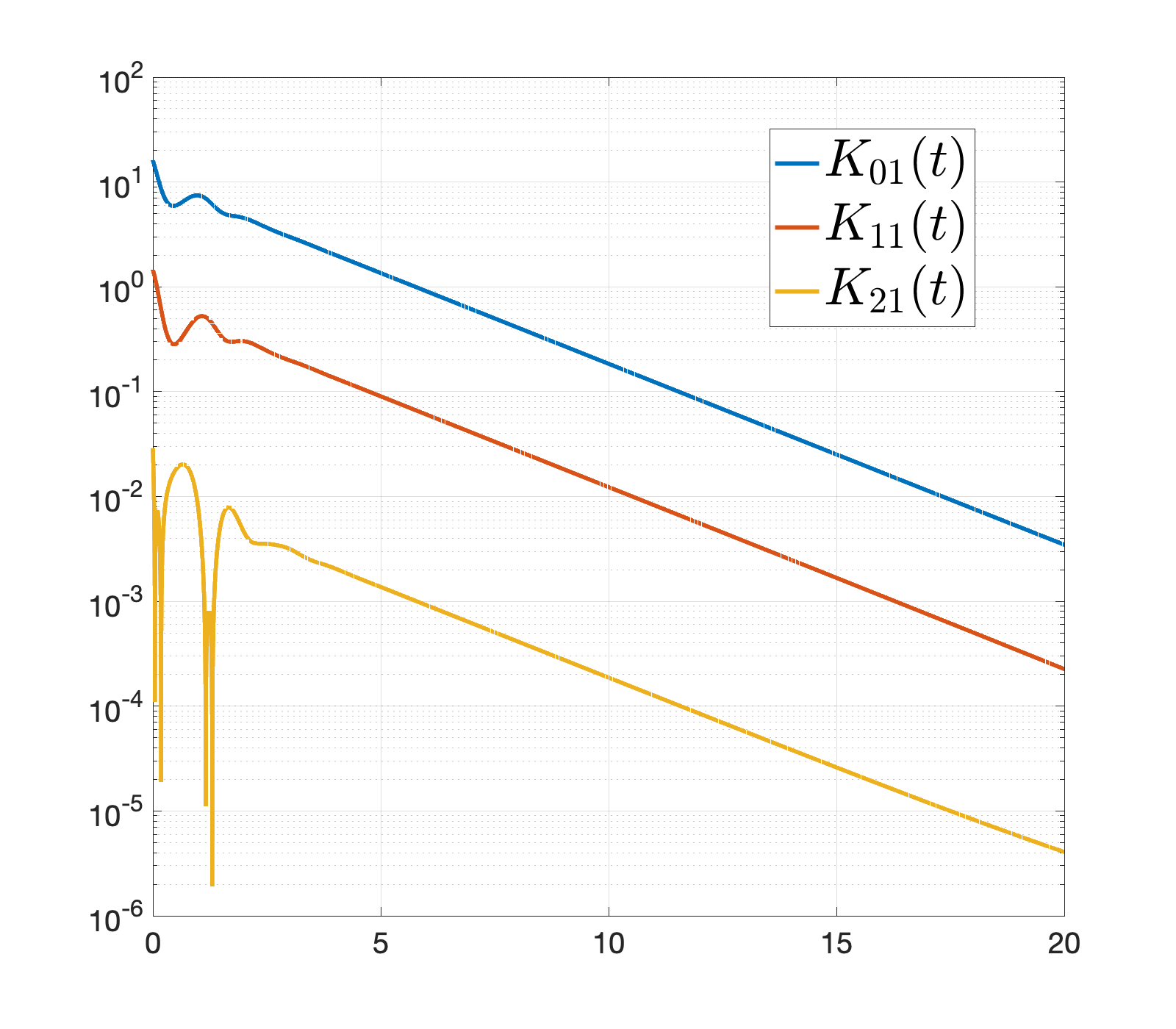}\label{HeQuK2}
}
\subfigure[$d = 3$]
{\includegraphics[width=0.31\textwidth]{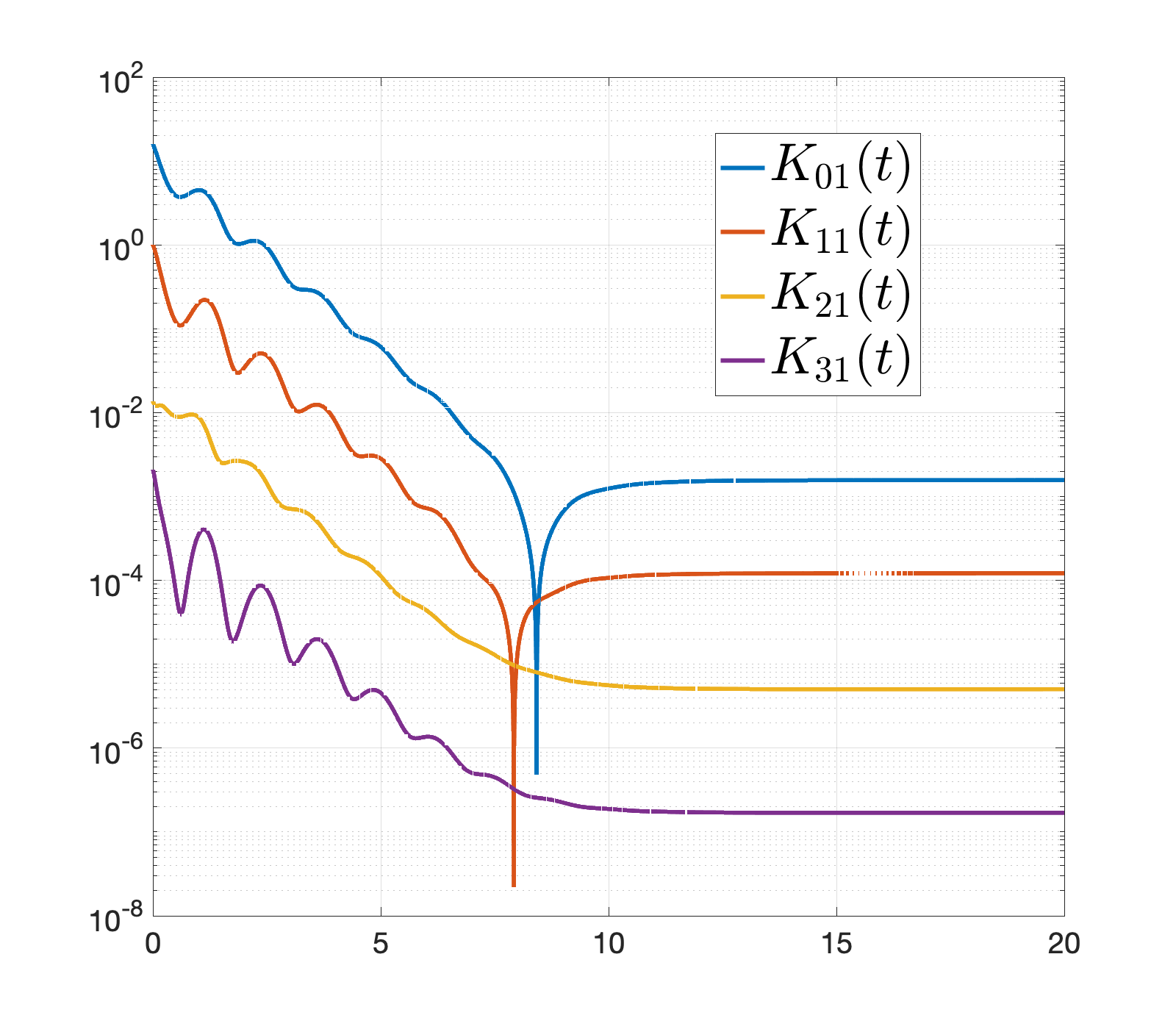}\label{HeQuK3}
}
\caption{Absolute values of the memory kernels $K_{l1}(t)$ for the nonlinear system \eqref{NLexample} for various values of $d$, the highest polynomial degree allowed in the expansion \eqref{MemKerExp}.}\label{HeQuK}
\end{figure}

For a simple demonstration of this approach, consider the two-dimensional nonlinear system
\eqn{
\phi_1'(t) &=& -\phi_1^2 + 8\phi_1\phi_2 \nonumber\\
\phi_2'(t) &=& \cos(\phi_1+\phi_2), \label{NLexample}
}
with initial conditions $\left(\phi_1(0),\phi_2(0)\right) = (x_1,x_2)$, and $\phi_1$ as the only resolved variable. Since
\eqn{
&& A_1(0 ; x_1,x_2) = 8x_1x_2 \nonumber\\
&\Rightarrow& 	\mathcal{L} A_1(0 ; x_1,x_2) = 8x_1\left[-x_1x_2 + 8x_2^2 + \cos(x_1+x_2)\right] \nonumber\\
&\Rightarrow& e^{t\mathcal{L}}\mathcal{L} A_1(0;x_1,x_2) = 8\phi_1\left[-\phi_1 \phi_2 + 8\phi_2^2 + \cos(\phi_1+\phi_2)\right] \label{NLExf}
}
and
\eqn{
[\mathcal{L}h_j](x_1,x_2) = \frac{1}{\sigma}H_j'\left(\frac{x_1 - \mu}{\sigma}\right) (-x_1^2 + 8x_1x_2) \Rightarrow e^{t\mathcal{L}}\mathcal{L}h_j (x_1,x_2) = \frac{1}{\sigma}H_j'\left(\frac{\phi_1 - \mu}{\sigma}\right) (-\phi_1^2 + 8\phi_1 \phi_2), \label{NLExg}
}
we can calculate $\left\{f_{l1} \right\}$ and $\left\{g_{lj}\right\}$ for $0 \leq l,j \leq d$ by solving the full system numerically for initial conditions of the form $\phi_1(0) = \mu + \sigma z_i$ and $\phi_2(0) = 0$. Figures \ref{HeQu1}, \ref{HeQu2} and \ref{HeQu3} show the results for the initial conditions considered in the previous section and for various values of $d$. In each case, we choose $\mu$ and $\sigma$ so the quadrature points are located in $[1,4]$, and employ a sufficiently high-order grid to ensure (empirical) convergence; we end up using $N_q = 20$, 30, and 40 for $d = 1$, 2, and 3 respectively. In Figure \eqref{HeQuK}, we also show the memory kernels $K_{l1}(t)$ computed for the various values of $d$. We note that these decay exponentially in time before plateauing off. 

We find that the memory kernels outperform the rudimentary Markovian system and that the performance improves when we retain more terms in the memory expansion. In particular, note the ability of the larger systems to correctly reproduce the limiting values as well as the dynamics before the steady state is achieved. This suggests that the inclusion of memory terms can help achieve significant gains in accuracy.

\section{Returning to viscous Burgers'}\label{SecAltVB}
The viscous Burgers' equation \eqref{ViscBurg} can be solved numerically by expanding the solution in terms of Fourier modes
\eqn{
u(t,z) = \sum_{k \in \mathbb{Z}} \theta_k(t) e^{ikz}, \quad u(0,z) = u_0(z), \nonumber
}
and solving the resulting system of ODEs
\eqn{
\theta_k'(t) = -\frac{ik}{2}\sum_{\substack{p,q \in \mathbb{Z} \\ p + q = k}} \theta_p(t) \theta_q(t) - \nu k^2\theta_k(t) =: R_k(\vec{\theta}), \label{VBODEs}
}
with $\theta_k(0) = \widehat{\left(u_0\right)}_k = \frac{1}{2\pi}\int_0^{2\pi} e^{-ikz}u_0(z) \ dz$, and where $\theta_{-k}(t) = \overline{\theta_k}(t)$ to ensure reality of $u$. In order to apply our machinery to this, we rewrite this system of ODEs in terms of the real and imaginary parts of $\theta_k(t)$, to wit, $\theta_k(t) = \phi_k(t) + i \psi_k(t)$, and obtain
\eqn{
\phi_k'(t) &=& k\sum_{\substack{p,q \in \mathbb{Z} \\ p + q = k}} \phi_p(t) \psi_q(t) - \nu k^2 \phi_k(t) =: T_k(\vec{\phi},\vec{\psi}) \label{VBrealPhiODEs}\\
\psi_k'(t) &=& -\frac{k}{2}\sum_{\substack{p,q \in \mathbb{Z} \\ p + q = k}} \left(\phi_p(t) \phi_q(t) -\psi_p(t) \psi_q(t)\right)- \nu k^2 \psi_k(t) =: S_k(\vec{\phi},\vec{\psi}), \label{VBrealPsiODEs}
}
with $\phi_{-k}(t) = \phi_k(t)$ and $\psi_{-k}(t) = -\psi_k(t)$ for all $k$. We also note that $\phi_0'(t) \equiv \psi_0'(t) \equiv 0$ so we have $\phi_0(t) \equiv \widehat{\left(u_0\right)}_0$ and $\psi_0(t) \equiv 0$. Since these modes do not evolve at all, we can omit these from our reduced-order model; more precisely, observe that the solution to \eqref{ViscBurg} can be obtained by solving
\eqn{
w_t = \pz\left(-\frac{1}{2} w^2\right) + \nu w_{zz}, \quad w(0,z) = w_0(z) := u_0(z) - \widehat{\left(u_0\right)}_0, \label{ViscBurgRed}
}
and setting $u(t,z) = \widehat{\left(u_0\right)}_0 + w\left(t,z - t\widehat{\left(u_0\right)}_0\right)$. Since
\eqn{
\widehat{\left(w_0\right)}_0 = \frac{1}{2\pi} \int_0^{2\pi} w_0(z) \ dz = \frac{1}{2\pi} \int_0^{2\pi}  u_0(z) - \widehat{\left(u_0\right)}_0 \ dz = 0, \nonumber
}
we can, without loss of generality, restrict ourselves to problems with a zero constant mode. In tandem with the reality constraints introduced earlier, this allows us to only focus only on the positively indexed modes. 

For a given $M > 0$, we take the resolved variables to be $\{\phi_{k}, \psi_{k}\}_{1 \leq k \leq M}$, yielding a total of $N_\text{rom} = 2M$ resolved modes. The resulting reduced-order models are
\eqn{
\phi_k'(t) &=&  \widecheck{T}_k(\widehat{\phi}(t),\widehat{\psi}(t)) + \sum_{j = 1}^{J} \int_0^t K^{(1)}_{jk}(t-s) h_j\left(\widehat{\phi}(s),\widehat{\psi}(s)\right) \ ds, \quad 1 \leq k \leq M, \nonumber\\
\psi_k'(t) &=& \widecheck{S}_k(\widehat{\phi}(t),\widehat{\psi}(t)) + \sum_{j = 1}^{J} \int_0^t K^{(2)}_{jk}(t-s) h_j\left(\widehat{\phi}(s),\widehat{\psi}(s)\right) \ ds, \quad 1 \leq k \leq M, \label{VBROMs}
}
where
\eqn{
\widecheck{T}_k(\widehat{x},\widehat{y}) &=&  k\sum_{\substack{|p|,|q| \leq M \\ p + q = k}} x_p y_q - \nu k^2 x_k , \qquad \widecheck{S}_k(\widehat{x},\widehat{y}) =  -\frac{k}{2}\sum_{\substack{|p|,|q| \leq M \\ p + q = k}} \left(x_px_q -  y_p y_q\right) - \nu k^2 y_k, \label{VBMarkODEs}
}
are the Markovian terms. We emphasize that we take $x_0 = y_0 = 0$ and $x_{-p} = x_p$ and $y_{-p} = -y_p$ for $p \geq 1$.

Our basis $\{h_j(\widehat{x},\widehat{y})\}$ over the resolved variables consists of up to cubic Hermite polynomials. Denoting the one-dimensional Hermite quadrature rule by $\{(r_i,\omega_i)\}_{1 \leq i \leq N_q}$, we solve the full system \eqref{VBrealPhiODEs} and \eqref{VBrealPsiODEs} for initial conditions chosen over all combinations of $\phi_k(0)/\sigma_k , \psi_l(0)/\sigma_l \in \{r_i\}_{1 \leq i \leq N_q}$ for $1 \leq k,l \leq M$, and zero otherwise. The $g_{lj}(t)$ terms in \eqref{fgdefs} can then be calculated by making use of
\eqn{
\left[\mathcal{L} h_j \right]\left(\widehat{x},\widehat{y}\right) = \sum_{k = 1}^M \frac{1}{\sigma_k}\left(\partial_{x_k}h_j \right)\left(\widehat{x},\widehat{y}\right) T_k(\vec{x},\vec{y}) + \sum_{l = 1}^M \frac{1}{\sigma_l} \left(\partial_{y_l}h_j \right)\left(\widehat{x},\widehat{y}\right) S_l(\vec{x},\vec{y}) \label{LiovH_j}
}
and applying the evolution operator $e^{t\mathcal{L}}$. For the $f_{lk}(t)$ terms, we note that the initial values of the noise corresponding to the resolved ``$\phi_k$'' are
\eqn{
A_k^{(1)}(0;\vec{x},\vec{y}) = T_k(\vec{x},\vec{y}) - \widecheck{T}_k(\widehat{x},\widehat{y}) = k\sum_{\substack{|p| \text{ or } |q| > M \\ p + q = k}} x_p y_q \nonumber
}
so that
\eqn{
\mathcal{L}A_k^{(1)}(0;\vec{x},\vec{y}) = k\sum_{\substack{|p| \text{ or } |q| > M \\ p + q = k}} y_qT_p(\vec{x},\vec{y}) +  x_qS_p(\vec{x},\vec{y}). \nonumber
}

Setting $\vec{v} = \vec{x} + i\vec{y}$ and noting that $\overrightarrow{R}(\vec{v}) = \overrightarrow{T}(\vec{x},\vec{y}) + i\overrightarrow{S}(\vec{x},\vec{y})$, we can write
\eqn{
y_qT_p(\vec{x},\vec{y}) +  x_qS_p(\vec{x},\vec{y}) = \text{Im}\left(v_qR_p(\vec{v})\right) \nonumber
}
with the result that
\eqn{
k\sum_{\substack{|p| \text{ or } |q| > M \\ p + q = k}} y_qT_p(\vec{x},\vec{y}) +  x_qS_p(\vec{x},\vec{y}) = -\text{Re}\left\{ik\sum_{\substack{|p| \text{ or } |q| > M \\ p + q = k}} v_qR_p(\vec{v}) \right\}, \nonumber
}
or alternatively
\eqn{
\mathcal{L}A_k^{(1)}(0;\vec{x},\vec{y}) = -\text{Re}\left\{ik\sum_{\substack{p,q \in \mathbb{Z} \\ p + q = k}} v_qR_p(\vec{v}) - ik\sum_{\substack{|p|,|q| \leq M \\ p + q = k}} v_qR_p(\vec{v})\right\}. \nonumber
}

The convolutional sums in this expression are more easily computed in the physical space; upon applying the evolution operator $e^{t\mathcal{L}}$, we obtain
\eqn{
e^{t\mathcal{L}}\mathcal{L}A_k^{(1)}(0;\vec{x},\vec{y}) = -\text{Re}\left\{\partial_z\left[u(-uu_z - \nu u_{zz})\right] - \partial_z\left[\left(\mathcal{P}_{M}u\right) \mathcal{P}_{M} (-uu_z - \nu u_{zz})\right]\right\}^{\widehat{}}_k. \label{VBf1term}
}

Similarly, for the noise corresponding to the resolved ``$\psi_k$'', we have $A^{(2)}_k(0;\vec{x},\vec{y}) = S_k(\vec{x},\vec{y}) - \widecheck{S}_k(\widehat{x},\widehat{y})$ and
\eqn{
e^{t\mathcal{L}}\mathcal{L}A_k^{(2)}(0;\vec{x},\vec{y}) = -\text{Im}\left\{\partial_z\left[u(-uu_z - \nu u_{zz})\right] - \partial_z\left[\left(\mathcal{P}_{M}u\right) \mathcal{P}_{M} (-uu_z - \nu u_{zz})\right]\right\}^{\widehat{}}_k. \label{VBf2term}
}

These can then be employed as described before to calculate the memory kernels $\{K_{jk}^{(1)}(t), K_{jk}^{(2)}(t)\}$ for $1 \leq k \leq M$ and $1 \leq j \leq J$, and which in turn can be used to numerically integrate the ROM \eqref{VBROMs}. 

To test this formulation, we calculate the memory kernels for the viscous Burgers' equation for $M = 3$ (so $N_\text{rom} = 6$) with $d = 1$ (``linear'') and $d = 3$ (``cubic'') (so $J = {7 \choose 1} = 7$ and $J = {9 \choose 3} = 84$ respectively). We take $\nu = 0.1$ and solve the full-order system up to $T = 20$ with a timestep of $10^{-4}$ using RK4. We further set $\sigma_k = e^{-k}$ and use $N_q = 4$ quadrature points for each resolved variable; thus, we need to calculate $(N_q)^{N_\text{rom}} = 4096$ full-order trajectories to calculate the $f_{lk}(t)$ and $g_{lj}(t)$ terms. The memory kernels (and the resulting ROMs) are solved for on the same time interval with a timestep of $10^{-3}$. 

\begin{figure}[tbph]
\centering
\subfigure[$u_0(z) = \sin(z)$]
{\includegraphics[width=0.31\textwidth]{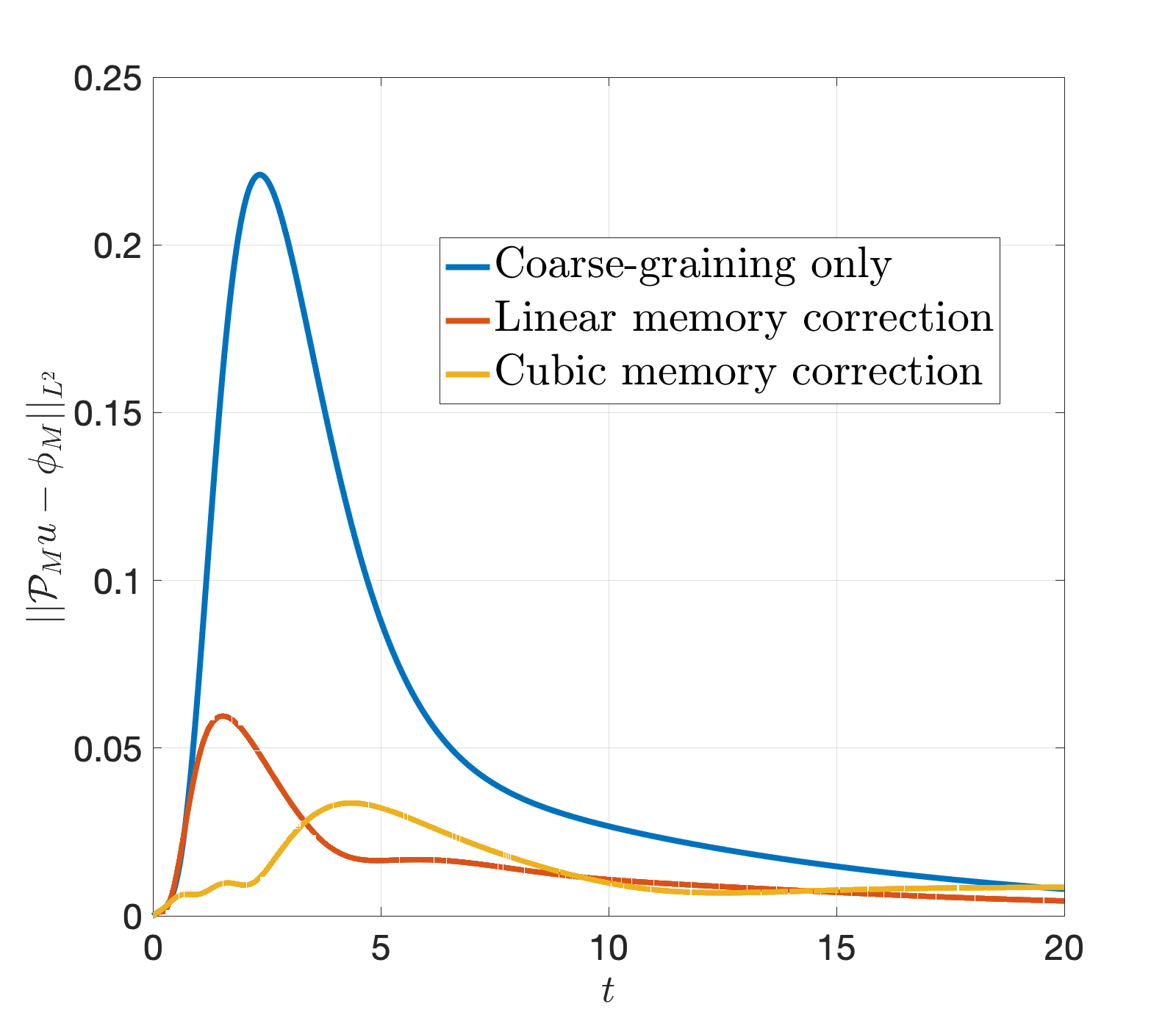}\label{MKCorr_sin}
}
\subfigure[$u_0(z) = \mathcal{P}_{M} \left(e^{\sin(z)}\right) $]
{\includegraphics[width=0.31\textwidth]{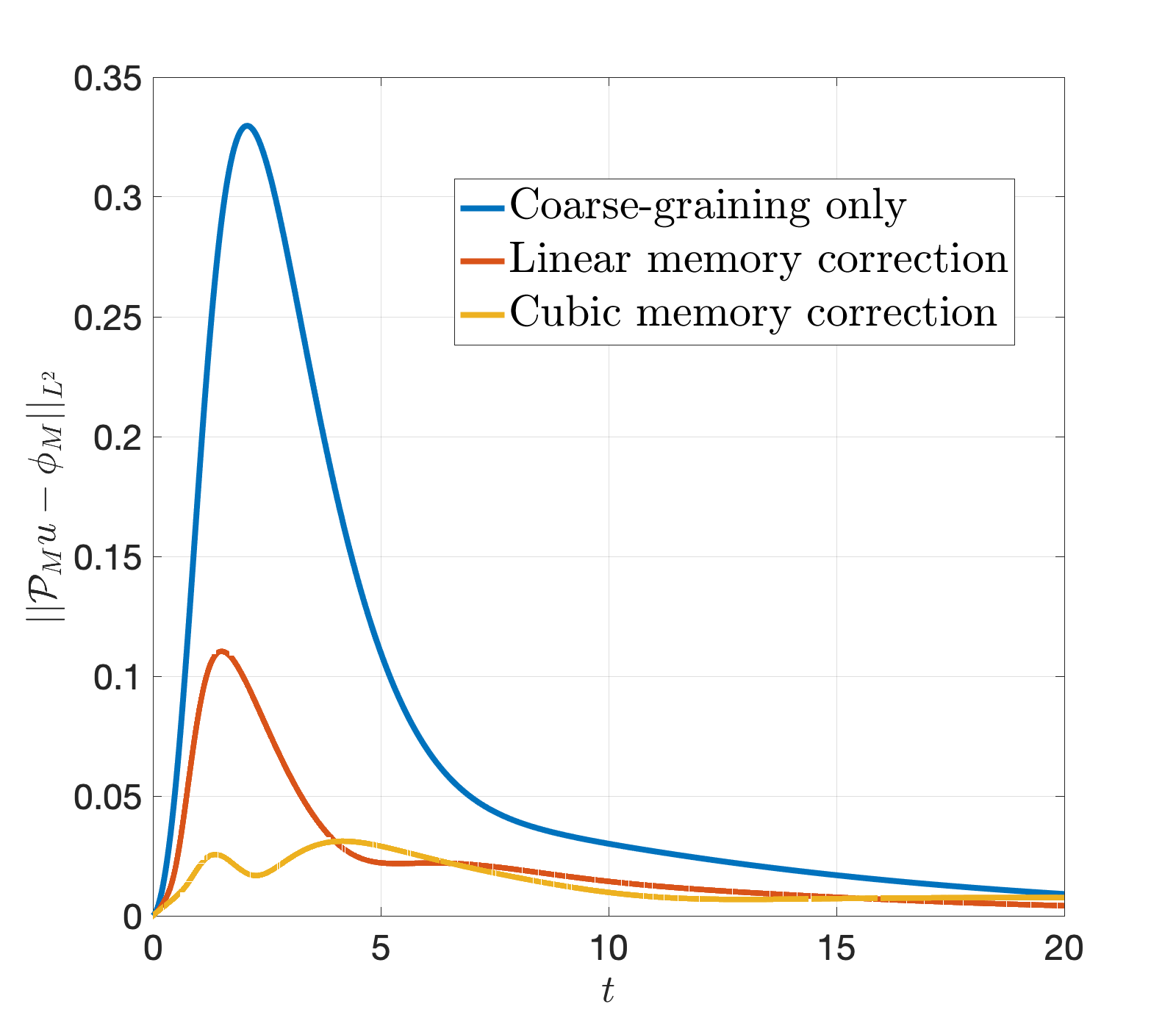}\label{MKCorr_expsin}
}
\subfigure[$u_0(z) = \mathcal{P}_{M} \left(\cos\left(2\sin(z)\right)\right) $]
{\includegraphics[width=0.31\textwidth]{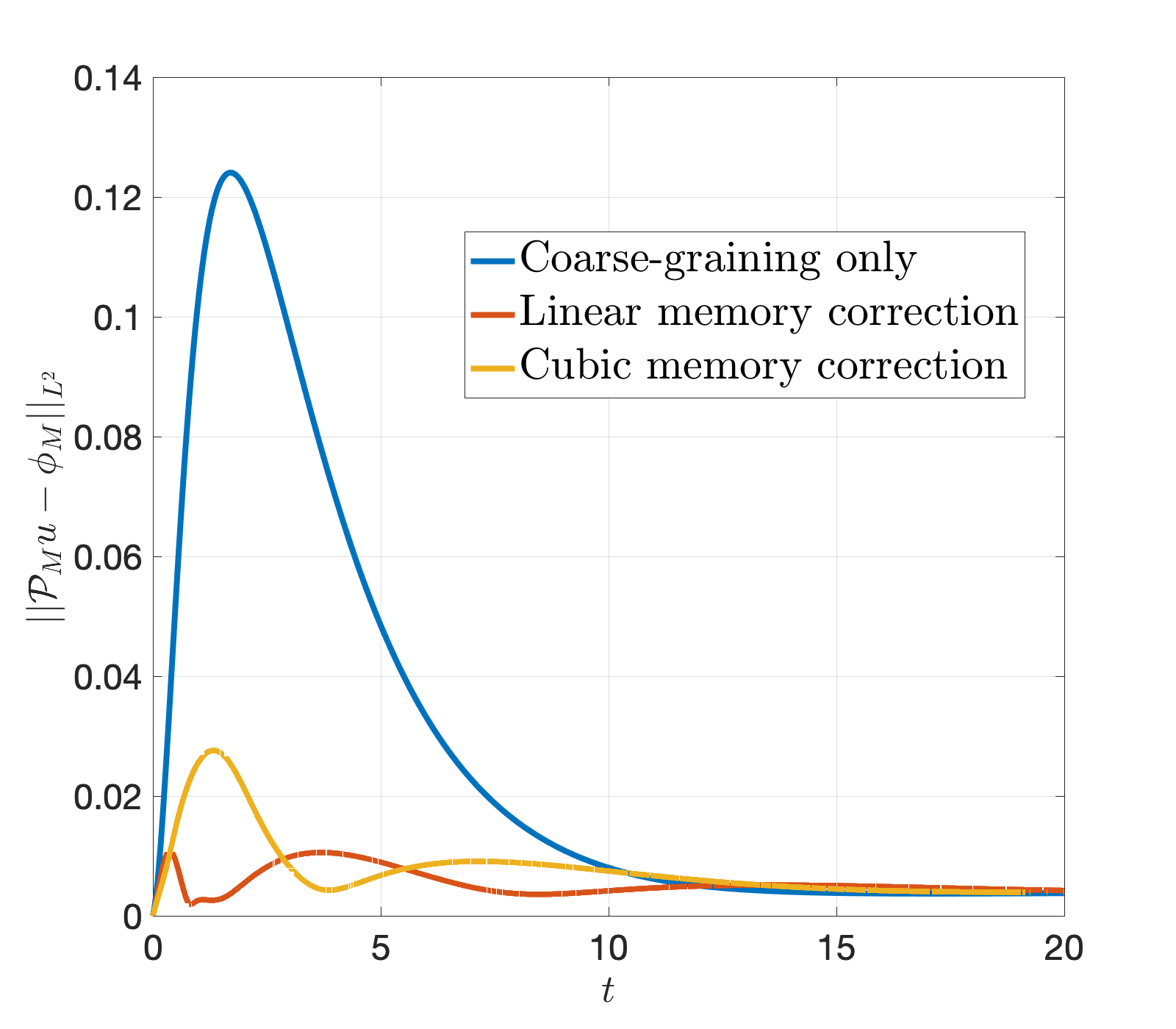}\label{MKCorr_cos2sin}
}
\caption{Error evolution $\norm{\mathcal{P}_{M}u_\text{ex} - u}_{L^2}$ for the solutions computed with and without the inclusion of the linear and cubic memory corrections to the stabilized Markovian system \eqref{MLSurrLP} for $M = 3$. We compare the solutions against the exact solution projected on the resolved modes as that represents (by definition) the best performance in $L^2$ we can hope to achieve using the allowed modes.}\label{MKCorr}
\end{figure}

In Figure \ref{MKCorr}, we display the results of adding the memory terms in the manner prescribed by \eqref{VBROMs} to the stabilized coarse-grained system \eqref{MLSurrLP} introduced in Section \ref{SecStabilizing}. The errors are computed by comparing the solutions of the Markovian and memory kernels solutions (both computed using RK2 with a timestep of $10^{-3}$) against $\mathcal{P}_{M}u_\text{ex}$, the exact solution projected on the resolved modes as that represents, by definition, the best performance in $L^2$ we can hope to achieve using only the allowed modes. We find that the incorporation of both linear and cubic memory corrections aids in improving the accuracy of the solutions across initial conditions.

\section{Discussion}\label{SecDisc}

Our preliminary results indicate that the instabilities in PDE--ML coupled systems can be avoided by coarse-graining the ``PDE modules'', i.e., parts of the system that are computed using conventional methods. This ensures that the inability of the ML surrogates to accurately mimic the high frequency components does not lead to a build-up of energy that cannot be suitably balanced or drained away. The resulting systems are then stable but possess only limited accuracy. We show that the use of memory-based correction terms stemming from the MZ formalism can aid in achieving better accuracy. Our future work aims to improve the efficiency of this approach, thereby permitting the use of larger values of $M$, and extending it to more complex problems. 

\bibliography{refs}

\end{document}